\RequirePackage{etoolbox}
\csdef{input@path}{{style/}{graphics/}}
\documentclass[seceqn,MSNbibl,number,citesort,dvips]{arxbj}
\usepackage{dsfont,mathbh,mathrsfs}
\usepackage{graphicx}

\aid{0}
\volume{21}
\issue{1}
\pubyear{2015}
\firstpage{1}
\lastpage{37}
\doi{10.3150/13-BEJ559} 

\makeatletter

\newcommand{\rright}{\right}
\newcommand{\lleft}{\left}

\renewcommand{\P}{\mathbb{P}}
\newcommand{\E}{\mathbb{E}}
\newcommand{\ddr}{\mathrm{d}}
\newcommand{\X}{\mathbb{X}}

\newproclaim{defi}{Definition}[section]
\newtheorem{lem}{Lemma}[section]
\newtheorem{cor}{Corollary}[section]
\newtheorem{prp}{Proposition}
\newremark{exe}{Example}[section]

\newcommand{\simiid}{\stackrel{\mathrm{i.i.d.}}{\sim}}

\newcommand{\Xcr}{\mathscr{X}}
\newcommand{\eqref}[1]{(\ref{#1})}
\makeatother

\begin{document}
\begin{frontmatter}

\title{Looking-backward probabilities for Gibbs-type exchangeable
random partitions}

\runtitle{Looking-backward probabilities}

\begin{aug}
\author[1]{\inits{S.}\fnms{Sergio} \snm{Bacallado}\thanksref{1}\ead[label=e1]{sergio.bacallado@gmail.com}},
\author[2,3]{\inits{S.}\fnms{Stefano} \snm{Favaro}\thanksref{2,3}\ead[label=e2]{stefano.favaro@unito.it}}
\and
\author[4]{\inits{L.}\fnms{Lorenzo} \snm{Trippa}\thanksref{4}\ead[label=e3]{ltrippa@jimmy.harvard.edu}}
\address[1]{Department of Statistics, Stanford University, Sequoia
Hall, Stanford, CA 94305, USA.\\ \printead{e1}}
\address[2]{Department of Economics and Statistics, University of
Torino, Corso Unione Sovietica 218/bis, 10134 Torino, Italy. \printead{e2}}
\address[4]{Harvard School of Public Health and Dana-Faber Cancer
Institute, 450 Brookline Avenue CLSB 11039 Boston, MA 02215, USA.
\printead{e3}}
\address[3]{Collegio Carlo Alberto, Via Real Collegio 30, 10024
Moncalieri, Italy}
\end{aug}

\received{\smonth{10} \syear{2012}}
\revised{\smonth{6} \syear{2013}}

%
\begin{abstract}
Gibbs-type random probability measures and the exchangeable random
partitions they induce
represent the subject of a rich and active literature. They provide a
probabilistic framework
for a wide range of theoretical and applied problems that are typically
referred to as species
sampling problems. In this paper, we consider the class of
looking-backward species sampling
problems introduced in Lijoi \textit{et al.} (\textit{Ann. Appl.
Probab.} \textbf{18} (2008)
1519--1547) in Bayesian nonparametrics. Specifically, given some
information on the random
partition induced by an initial sample from a Gibbs-type random
probability measure, we study
the conditional distributions of statistics related to the old species,
namely those species detected
in the initial sample and possibly re-observed in an additional sample.
The proposed results contribute
to the analysis of conditional properties of Gibbs-type exchangeable
random partitions, so far focused
mainly on statistics related to those species generated by the
additional sample and not already detected in the initial sample.
\end{abstract}

%

\end{frontmatter}

\section{Introduction}

Let $\mathbb{X}$ be a complete and separable metric space equipped
with the Borel $\sigma$-algebra $\mathscr X$, and let $(X_{i})_{i\geq
1}$ be an exchangeable sequence of $\mathbb{X}$-valued random
variables defined on some probability space $(\Omega,\mathscr
{F},\mathds{P})$. According to the celebrated de Finetti's
representation theorem there exists a random probability measure
$\tilde{P}$ on $\mathbb{X}$ such that, conditionally on $\tilde{P}$,
the random variables $(X_{i})_{i\geq1}$ are independent and
identically distributed according to $\tilde{P}$, that is,
\begin{eqnarray*}
X_i | \tilde P & \simiid& \tilde P,
\\
\tilde P & \sim& \Pi.
\end{eqnarray*}
The distribution $\Pi$ is commonly known as the de Finetti probability
measure of $(X_{i})_{i\geq1}$ and it takes on the interpretation of
the prior distribution in Bayesian nonparametrics. In the present
paper, we consider almost surely discrete random probability measures,
namely $\tilde P$ is such that $\Pi[\tilde P\in\mathscr{D}]=1$,
where $\mathscr{D}$ stands for the set of discrete probability
measures on $(\X,\Xcr)$.

If $\tilde{P}$ is discrete almost surely, we expect ties in a sample
$(X_{1},\ldots,X_{n})$ from $\tilde{P}$; that is, we expect
$K_{n}\leq n$ distinct observations with frequencies $\mathbf
{N}_{n}=(N_{1},\ldots,N_{K_{n}})$ satisfying $\sum_{1\leq i\leq
K_{n}}N_{i}=n$. Accordingly, the sample induces a random partition of
the set $\{1,\ldots,n\}$, in the sense that any index $i\neq j$
belongs to the same partition set if and only if $X_{i}=X_{j}$. We
denote by $p_{j}^{(n)}(n_{1},\ldots,n_{j})$ the symmetric function
corresponding to the probability of any particular partition of $\{
1,\ldots,n\}$ having $j$ distinct blocks with frequencies
$(n_{1},\ldots,n_{j})$. This function is known as the exchangeable
partition probability function (EPPF), a concept introduced in \cite
{Pit95} as a development of earlier results in \cite{Kin78}. The
EPPF can be specified for every $n\geq1$ and $1\leq j\leq n$ 
either via the exchangeable sequence $(X_{i})_{i\geq1}$ or by defining
a random partition of $\mathbb{N}$. In the latter case, the
distribution of the random partition must satisfy certain consistency
conditions and a symmetry property that guarantees exchangeability. See
\cite{Pit06} and references therein for a comprehensive account on EPPFs.

Exchangeable random partitions play an important role in a variety of
research areas. In population genetics, models for exchangeable random
partitions are useful for describing the configurations of a sample of
genes into a number of distinct allelic types. See \cite{Ewe98} and
references therein. In machine learning, probabilistic models for
linguistic applications are often based on clustering structures for
collections of words in documents. See, for example, \cite{Teh06}
and \cite{bnp_teh} for a review. In Bayesian nonparametrics,
exchangeable random partitions are commonly employed at the latent
level of complex hierarchical mixture models. See \cite{bnp_lij} and
references therein for a review. Other areas of application include
storage problems, excursion theory, combinatorics, number theory and
statistical physics. Broadly speaking, exchangeable random partitions
and their associated EPPFs provide a flexible probabilistic framework
for a wide range of theoretical and applied problems that are typically
referred to as species sampling problems, namely problems concerning a
population composed of individuals belonging to different species.
Indeed, the number of partition blocks $K_{n}$ take on the
interpretation of the number of distinct species in the sample
$(X_{1},\ldots,X_{n})$ and the $N_{i}$'s are the corresponding species
frequencies. Given the relevance and intuitiveness of such a framework,
throughout the paper we will resort to the species metaphor.

The main object of our investigation is the class of Gibbs-type
exchangeable random partitions. These are random partitions which arise
by sampling from a random probability measure, say of Gibbs-type, here
denoted by $\tilde{P}_{G}$. See \cite{Pit03} for details.
Introduced in \cite{Gne05} these exchangeable random partitions
represent the subject of a rich and active literature. A recent
development, first proposed in \cite{Lij308}, is the study of their
conditional properties. This study consists in evaluating, conditional
on some information about the random partition induced by an initial
sample $(X_{1},\ldots,X_{n})$ from $\tilde{P}_{G}$, the distribution
of certain statistics of an additional sample $(X_{n+1},\ldots
,X_{n+m})$. In particular, in \cite{Lij308} the main focus is on
the conditional distributions of statistics related to the \textit
{new} species, namely those species generated by the additional sample
and not coinciding with species already detected in the initial sample.
A representative example is given by the distribution of the number of
new distinct species generated by $(X_{n+1},\ldots,X_{n+m})$,
conditional on the information of both the number of distinct species
in $(X_{1},\ldots,X_{n})$ and their corresponding frequencies. See
\cite{Fav12a} for a generalization to the number of new distinct\vadjust{\goodbreak}
species with a certain frequency of interest. As shown in \cite
{Fav12a,Lij207} and \cite{Lij308} these conditional
distributions have direct applications in Bayesian nonparametric
analysis of species sampling problems arising in ecology and genomics.
We refer to \cite{Deb13a,Deb13b,Fav12b} and
\cite{gri07} for other contributions at the interface between
Bayesian nonparametrics and Gibbs-type exchangeable random partitions.

Many problems in the conditional analysis of Gibbs-type exchangeable
random partitions remain unresolved. For instance, \cite{Lij308}
pointed out the practical interest in the conditional distributions of
statistics related to the \textit{old} species, namely those species
detected in the initial sample and possibly re-observed in the
additional sample. Two illustrative examples are given in Proposition~4
of \cite{Lij308} and in Theorem~3 of \cite{Fav12a}. In general
the class of species sampling problems concerning old species has been
referred to as \textit{looking-backward} and it represent the focus of
the present paper. We study two novel, and practically applicable,
looking-backward species sampling problems. In particular, we derive
\begin{enumerate}[(ii)]
\item[(i)] the conditional distribution of the number of old distinct
species re-observed in $(X_{n+1},\ldots,\allowbreak X_{n+m})$, given complete or
incomplete information on the random partition induced by
$(X_{1},\ldots,X_{n})$;
\item[(ii)] the conditional distribution of the number of old distinct
species re-observed with a specific frequency of interest in
$(X_{n+1},\ldots,X_{n+m})$, given complete or incomplete information
on the random partition induced by $(X_{1},\ldots,X_{n})$.
\end{enumerate}
Specifically, by complete information we refer jointly to the number of
distinct species in $(X_{1},\ldots,X_{n})$ and their frequencies,
whereas by incomplete information we refer solely to the number of
distinct species in $(X_{1},\ldots,X_{n})$. Besides the sets of
complete and incomplete information, we also consider almost-complete
information. This information refers jointly to the number of distinct
species in $(X_{1},\ldots,X_{n})$ and a subset of their corresponding
frequencies.

The present paper broadens the scope of previous literature on
conditional distributions for Gibbs-type exchangeable random
partitions, by investigating in depth some statistics related to old
species. In the framework of Gibbs-type exchangeable random partitions,
looking-backward problems create a distinction between conditioning on
complete, incomplete and almost complete information, which to the best
of our knowledge has not been dealt with explicitly in previous
studies. We expect the results introduced here to have an impact in the
analysis of Bayesian nonparametric models for species sampling
problems, which have acquired increasingly complex forms in recent
years to meet the demands of scientific applications. The paper is
structured as follows. Section~\ref{sec2} recalls the definition of Gibbs-type
exchangeable random partition and introduces preliminary results
relevant to the analysis of their conditional structure. Section~\ref{sec3}
deals with the looking-backward species sampling problems (i) and (ii) in
the general case of Gibbs-type exchangeable random partitions and in
the special case of the celebrated Ewens--Pitman sampling model. The
context of almost-complete information is also dealt with in Section~\ref{sec3}.
Section~\ref{sec4} contains some numerical illustrations of the present results.
Proofs are deferred to the \hyperref[app]{Appendix}.


\section{Preliminaries and main definitions}\label{sec2}
Gibbs-type exchangeable random partitions were introduced in \cite
{Gne05} and further investigated in \cite{Pit03}. This class of
exchangeable random partitions is characterized by an EPPF with a
product form, a feature which is crucial for mathematical tractability
and, in particular, facilitates intuition. Let $\mathcal{D}_{n,j}=\{
(n_{1},\ldots,n_{j}): n_{i}\geq1\mbox{ and }\sum_{i=1}^{j}n_{i}=n\}
$ be the set of the partitions of $n\geq1$ into $j\leq n$ positive
integers. Moreover, for any $x>0$ and any positive integer $n$, we
denote by $(x)_{n\uparrow1}$ and $(x)_{n\downarrow1}$ the rising
factorials and falling factorials, respectively.
%
\begin{defi}\label{defi_gibbs}
Let $(X_{i})_{i\geq1}$ be an exchangeable sequence directed by $\tilde
{P}_{G}$. Then, the exchangeable random partition induced by
$(X_{i})_{i\geq1}$ is said of Gibbs-type and it is characterized by an
EPPF of the form
%
\begin{equation}
\label{eq:gibbs_defi} p_{j}^{(n)}(n_{1},
\ldots,n_{j})=V_{n,j}\prod_{i=1}^{j}(1-
\sigma )_{(n_{i}-1)\uparrow1},
\end{equation}
for $\sigma<1$ and nonnegative weights $(V_{n,j})_{j\leq n,n\geq1}$
satisfying the recursion $V_{n,j}=V_{n+1,j+1}+(n-j\sigma)V_{n+1,j}$,
with $V_{1,1}=1$.
\end{defi}

Gibbs-type exchangeable random partitions are completely specified by
the parameter $\sigma<1$ and the collection of weights
$(V_{n,j})_{j\leq n,n\geq1}$ satisfying a backward recursion. Note
that Definition~\ref{defi_gibbs} also provides the distribution of the
number $K_{n}$ of distinct species in a sample of size $n$ from $\tilde
{P}_{G}$, that is,
%
\begin{equation}
\label{eq:gibbs_dist} \P[K_{n}=j]=V_{n,j}\frac{\mathscr{C}(n,j;\sigma)}{\sigma^{j}},
\end{equation}
with $\mathscr{C}(n,j;\sigma)$ being the so-called generalized
factorial coefficient. We refer to \cite{Cha05} for details. The
next example recalls the Ewens--Pitman sampling model, a noteworthy
example of Gibbs-type exchangeable random partition introduced in
\cite{Pit95} and generalizing the celebrated Ewens sampling model
in \cite{Ewe72}. See \cite{Arr03} and references therein for a
comprehensive account on the Ewens sampling model. Another notable
Gibbs-type exchangeable random partition, still related to the
Ewens--Pitman sampling model, has been recently introduced and
investigated in~\cite{Gne10}.
%
\begin{exe}
For any $\sigma\in(0,1)$ and $\theta>-\sigma$, the Ewens--Pitman
sampling model is a Gibbs-type exchangeable random partition with
weights $(V_{n,j})_{j\leq n,n\geq1}$ of the following form
%
\begin{equation}
\label{eq:defi_ewpit} V_{n,j}=\frac{\prod_{i=0}^{j-1}(\theta+i\sigma)}{(\theta
)_{n\uparrow1}}.
\end{equation}
The Ewens sampling model with parameter $\theta>0$ is recovered from
the Ewens--Pitman sampling model by letting $\sigma\rightarrow0$. See,
for example, \cite{Pit95} and \cite{Pit97} for details and
further developments.
\end{exe}

The recursion in Definition~\ref{defi_gibbs}, for a fixed $\sigma$,
cannot be solved in a unique way. The solutions form a convex set where
each element is the distribution of an exchangeable random partition.
Theorem~12 in \cite{Gne05} describes the extreme points of such a
convex set. For any $n\ge1$ let
%
\begin{eqnarray*}
c_{n}(\sigma)=\lleft\{ %
\begin{array} {l@{\qquad}l} 1& \mbox{if }
\sigma\in(-\infty,0),
\\
\log(n)& \mbox{if }\sigma=0,
\\
n^{\sigma}& \mbox{if }\sigma\in(0,1). \end{array} %
\rright.
\end{eqnarray*}
For every Gibbs-type exchangeable random partition there exists a
positive and almost surely finite random variable $S_{\sigma}$ such that
%
\begin{eqnarray*}
\frac{K_{n}}{c_{n}(\sigma)}\stackrel{\mathrm{a.s.}} {\longrightarrow }S_{\sigma},
\end{eqnarray*}
as $n\rightarrow+\infty$, and such that a Gibbs-type exchangeable
random partition is a unique mixture over
$\varkappa$ of extreme exchangeable random partitions for which
$S_{\sigma}=\varkappa$ almost surely.
For $\sigma\in(-\infty,0)$ the extremes are Ewens--Pitman sampling
models with parameter $(\sigma,-\sigma\varkappa)$; for $\sigma=0$
the extremes are Ewens sampling models with parameter $\varkappa\geq
0$; for $\sigma\in(0,1)$ the Ewens--Piman sampling models are not
extremes. See Section~6.1 in \cite{Pit03} for details on $S_{\sigma}$.

A generalization of Definition~\ref{defi_gibbs} has been recently
introduced in \cite{Lij308} to study conditional properties of
Gibbs-type exchangeable random partitions. To recall this
generalization a few quantities, analogous to those describing the
random partition induced by an initial sample $(X_{1},\ldots,X_{n})$
from $\tilde{P}_{G}$, need to be introduced. Let $X_{1}^{\ast},\ldots
,X_{K_{n}}^{\ast}$ be the labels identifying the $K_{n}$ distinct
species detected in the initial sample and, for any $m>1$, define
%
\begin{equation}
\label{eq:new_observ} L_{m}^{(n)}=\sum_{i=1}^{m}
\prod_{j=1}^{K_{n}}\mathbh{1}_{\{
X_{j}^{\ast}\}^{C}}(X_{n+i})
\end{equation}
as the number of observations in an additional sample $(X_{n+1},\ldots
,X_{n+m})$ not coinciding with any of the $K_{n}$ distinct species.
Denote by $K_{m}^{(n)}$ the number of new distinct species generated by
these $L_{m}^{(n)}$ observations and by $X_{K_{n}+1}^{\ast},\ldots
,X_{K_{n}+K_{m}^{(n)}}^{\ast}$ their corresponding identifying labels.
Therefore,
\begin{eqnarray*}
\mathbf{M}_{L_{m}^{(n)}}= (M_{1},\ldots,M_{K_{m}^{(n)}} ),
\end{eqnarray*}
with
%
\begin{equation}
\label{eq:freq_old} M_{i}=\sum_{j=1}^{m}
\mathbh{1}_{\{X_{K_{n}+i}^{\ast}\}}(X_{n+j})
\end{equation}
for $i=1,\ldots,K_{m}^{(n)}$, are the frequencies of the new
$K_{m}^{(n)}$ distinct species detected among the $L_{m}^{(n)}$
observations of the additional sample. Analogously,
%
\begin{eqnarray*}
\mathbf{S}_{m-L_{m}^{(n)}}= (S_{1},\ldots,S_{K_{n}} ),
\end{eqnarray*}
with
%
\begin{equation}
\label{eq:freq_new} S_{i}=\sum_{j=1}^{m}
\mathbh{1}_{\{X^{\ast}_{i}\}}(X_{n+j}),
\end{equation}
corresponds to number of observations, among the $m-L_{m}^{(n)}$
observations of the additional sample, coinciding with the $i$th
distinct old species detected in the initial sample, for $i=1,\ldots
,K_{n}$. As pointed out in \cite{Lij207}, from a Bayesian
nonparametric perspective the joint conditional distribution of the
random variables \eqref{eq:new_observ}, \eqref{eq:freq_old}, \eqref
{eq:freq_new} and $K_{m}^{(n)}$, given $(X_{1},\ldots,X_{n})$, can be
interpreted as the posterior counterpart of the EPPF \eqref
{eq:gibbs_defi}. This then provides a natural framework for Bayesian
nonparametric analysis of species sampling problem.

In \cite{Lij308}, the main focus is on conditional distributions of
statistics related to the new species generated by $(X_{n+1},\ldots
,X_{n+m})$. For instance, by suitably marginalizing the joint
conditional distribution of the random variables \eqref
{eq:new_observ}, \eqref{eq:freq_old}, \eqref{eq:freq_new} and
$K_{m}^{(n)}$, given $(X_{1},\ldots,X_{n})$, one obtains the
conditional distribution of the number of new distinct species, namely
%
\begin{equation}
\label{eq:distinct} \P \bigl[K_{m}^{(n)}=k | K_{n}=j,
\mathbf{N}_{n}=\mathbf{n} \bigr]=\frac{V_{n+m,j+k}}{V_{n,j}}\frac{\mathscr{C}(m,k;\sigma
,-n+j\sigma)}{\sigma^{k}}
\end{equation}
with $\mathscr{C}(n,j;\sigma,\rho)$ being the so-called noncentral
generalized factorial coefficient. We refer to \cite{Cha05} for
details. Accordingly, the Bayesian nonparametric estimator, under
quadratic loss function, of the number of new distinct species
generated by the additional sample coincides with
%
\begin{equation}
\label{eq:estimator_dist} \mathcal{K}_{m}^{(n)}=\E \bigl[K_{m}^{(n)}
| K_{n}=j,\mathbf {N}_{n}=\mathbf{n} \bigr]= \E
\bigl[K_{m}^{(n)} | K_{n}=j \bigr].
\end{equation}
We refer to \cite{Deb13b,Lij207,Lij07Bio}
and \cite{Lij308} for applications of \eqref{eq:distinct} and
\eqref{eq:estimator_dist}, under the choice of $V_{n,j}$ in \eqref
{eq:defi_ewpit}, to Bayesian nonparametric inference for species
variety in genetic experiments. As a generalization of \eqref
{eq:distinct}, Theorem~3 in \cite{Fav12a} provides the conditional
distribution, given $(X_{1},\ldots,X_{n})$, of
%
\begin{equation}
\label{eq:estimator_dist_freq} \sum_{i=1}^{K_{n}}
\mathbh{1}_{\{N_{i}+S_{i}=l\}}+\sum_{i=1}^{K_{m}^{(n)}}
\mathbh{1}_{\{M_{i}=l\}},
\end{equation}
for any $l=1,\ldots,n+m$. In words, \eqref{eq:estimator_dist_freq}
corresponds to the number of distinct species with frequency $l$
generated by $(X_{n+1},\ldots,X_{n+m})$. The conditional expected
value of \eqref{eq:estimator_dist_freq}, given $(X_{1},\ldots
,X_{n})$, provides the Bayesian nonparametric estimator, under
quadratic loss function, of the number of distinct species with
frequency $l$ generated by the additional sample.


\section{Two looking-backward probabilities}\label{sec3}

Before presenting our results, it is worth stating the fundamental
difference between looking-backward species sampling problems and the
species sampling problems\vadjust{\goodbreak} investigated in \cite{Lij308}. A~common
feature of the conditional distributions introduced in \cite
{Lij308} is their independence from the information on the
frequencies $\mathbf{N}_{n}$ induced by the initial sample
$(X_{1},\ldots,X_{n})$. As a representative example, note that the
distribution \eqref{eq:distinct} satisfies the following identity
\begin{eqnarray*}
\P \bigl[K_{m}^{(n)}=k | K_{n}=j,
\mathbf{N}_{n}=\mathbf{n} \bigr]=\P \bigl[K_{m}^{(n)}=k
| K_{n}=j \bigr].
\end{eqnarray*}
Such a property of independence characterizes all the statistics
concerning the new species in the additional sample $(X_{n+1},\ldots
,X_{n+m})$. Indeed, since \eqref{eq:freq_new} does not contain any
information on new species, the conditional distributions of these
statistics can be obtained from the joint conditional distribution of
the random variables \eqref{eq:new_observ}, \eqref{eq:freq_old} and
$K_{m}^{(n)}$, given $(X_{1},\ldots,X_{n})$. In Proposition~1 of
\cite{Lij207}, this joint conditional distribution is shown to be
independent of $\mathbf{N}_{n}$. Hence, $K_{n}$ is a sufficient
statistic for the species sampling problems discussed in~\cite{Lij308}.

Differently, the conditional distributions of statistics concerning old
species depend on the information of both the number $K_{n}$ of
distinct species and the corresponding frequencies $\mathbf{N}_{n}$.
This is to say that, letting $T_{m}^{(n)}$ be a statistic related to
old species, in most cases, one obtains
%
\begin{equation}
\label{eq:nosufficiency} \P \bigl[T_{m}^{(n)}\in\cdot|
K_{n}=j,\mathbf{N}_{n}=\mathbf {n} \bigr]\neq\P
\bigl[T_{m}^{(n)}\in\cdot| K_{n}=j \bigr].
\end{equation}
As an example, the distribution of \eqref{eq:estimator_dist_freq}
satisfies \eqref{eq:nosufficiency}. See Theorem~3 of \cite{Fav12a}
for details. See also Proposition~4 in \cite{Lij308} for another
example. According to \eqref{eq:nosufficiency}, the analysis of the
looking-backward species sampling problems naturally leads to consider
at least two sets of information on the random partition induced by
$(X_{1},\ldots,X_{n})$:\vadjust{\goodbreak} (i) a complete information, namely $K_{n}$ and
$\mathbf{N}_{n}$; (ii) an incomplete information, namely $K_{n}$. We
also consider almost-complete information, namely $K_{n}$ and a subset
of $\mathbf{N}_{n}$. In the next subsections, we present and discuss
the results of our paper. We focus on deriving the conditional
distributions of two looking-backward statistics, given complete or
incomplete information. This will be the subject of Section~\ref{subsec1} and Section~\ref{subsec2}. The conditional distributions
of these two statistics given almost-complete information can be
derived through similar arguments applied when conditioning on
incomplete information. This will be discussed in Section~\ref{subsec3}.

\subsection{Probabilities of re-observing old species}\label{subsec1}

In this section, we consider the distribution of the number of old
distinct species that are re-observed in $(X_{n+1},\ldots,X_{n+m})$,
conditional on complete and incomplete information on the random
partition induced by $(X_{1},\ldots,X_{n})$. Formally, in the context
of complete information, we are interested in the random variable
$R_{m}^{(n,j,\mathbf{n})}$ which is defined in distribution as
%
\begin{equation}
\label{eq:def_main_1} \P\bigl[R_{m}^{(n,j,\mathbf{n})}=x\bigr]=\P \Biggl[\sum
_{i=1}^{K_{n}}\mathbh {1}_{\{S_{i}>0\}}=x \Bigl|
K_{n}=j,\mathbf{N}_{n}=\mathbf{n} \Biggr].
\end{equation}
In the context of incomplete information, we are interested in the
random variable $\tilde{R}_{m}^{(n,j)}$ which is defined in
distribution as
%
\begin{equation}
\label{eq:def_main_1_part} \P\bigl[R_{m}^{(n,j)}=x\bigr]=\P \Biggl[\sum
_{i=1}^{K_{n}}\mathbh{1}_{\{
S_{i}>0\}}=x \Bigl|
K_{n}=j \Biggr].
\end{equation}
In the next theorem, we derive the factorial moments of the random
variables in \eqref{eq:def_main_1} and \eqref{eq:def_main_1_part}. By
means of Theorem~1 in \cite{Fav12a}, we obtain \eqref{teo1_eq1}.
Accordingly, \eqref{teo1_eq2} follows from \eqref{teo1_eq1} by
suitably marginalizing the frequencies $\mathbf{N}_{n}$. These moments
then lead to the corresponding distributions by means of standard
arguments involving probability generating functions.
%
\begin{thm}\label{teo1}
Let $(X_{i})_{i\geq1}$ be an exchangeable sequence directed by $\tilde
{P}_{G}$. Then, for any integer $r\geq1$ one has
%
\begin{eqnarray}
\label{teo1_eq1} &&\E \bigl[ \bigl(R_{m}^{(n,j,\mathbf{n})}
\bigr)_{r\downarrow
1} \bigr]
\nonumber\\
&&\quad
=r!\sum_{v=0}^{r}{j-v\choose
r-v}(-1)^{v}
\\
&&\qquad {}
\times\sum_{\{c_{1},\ldots,c_{v}\}\in\mathcal
{C}_{j,v}}\sum
_{k=0}^{m}\frac{V_{n+m,j+k}}{V_{n,j}}\frac{\mathscr
{C}(m,k;\sigma,-n+\sum_{i=1}^{v}n_{c_{i}}+(j-v)\sigma)}{\sigma^{k}}\nonumber
\end{eqnarray}
and
%
\begin{eqnarray}
\label{teo1_eq2} &&\E \bigl[ \bigl(R_{m}^{(n,j)}
\bigr)_{r\downarrow1} \bigr]
\nonumber\\
&&\quad
=\frac{r!}{\mathscr{C}(n,j;\sigma)} \sum_{v=0}^{r}{j-v
\choose r-v}(-1)^{v}\sum_{s=v}^{n-(j-v)}{n
\choose s}\mathscr{C}(s,v;\sigma)\mathscr{C}(n-s,j-v;\sigma)
\\
&&\qquad {}
\times\sum_{k=0}^{m}
\frac{V_{n+m,j+k}}{V_{n,j}}\frac
{\mathscr{C}(m,k;\sigma,-n+s+(j-v)\sigma)}{\sigma^{k}},\nonumber
\end{eqnarray}
where $\mathcal{C}_{j,v}$ denotes the set of the $v$-combinations
(without repetitions) of the elements $\{1,\ldots,j\}$.
\end{thm}

The distributions of $R_{m}^{(n,j,\mathbf{n})}$ and $R_{m}^{(n,j)}$
are interpretable as the posterior distributions of the number of old
distinct species that are re-observed in $(X_{n+1},\ldots,X_{n+m})$
given, respectively, complete and incomplete information on the random
partition induced by $(X_{1},\ldots,X_{n})$. Accordingly, the Bayesian
nonparametric estimators, under a quadratic loss function, coincide
with the expected values of the random variables $R_{m}^{(n,j,\mathbf
{n})}$ and $R_{m}^{(n,j)}$. An expression for these Bayesian
nonparametric estimators, denoted by $\mathcal{R}_{m}^{(n,j,\mathbf
{n})}=\E[R_{m}^{(n,j,\mathbf{n})}]$ and $\mathcal{R}_{m}^{(n,j)}=\E
[R_{m}^{(n,j)}]$, is presented in the next corollary. See Proposition~\ref{prop1_1} and Proposition~\ref{prop1_2} for an expression of
these estimators under the Ewens--Pitman sampling model.
%
\begin{cor}
The Bayesian nonparametric estimator of the number of old distinct
species that are re-observed in an additional sample of size $m$, given
complete information on $(X_{1},\ldots,X_{n})$, coincides with
\begin{eqnarray*}
\mathcal{R}_{m}^{(n,j,\mathbf{n})}=j-\sum_{i=1}^{n}m_{i}
\sum_{k=0}^{m}\frac{V_{n+m,j+k}}{V_{n,j}}
\frac{\mathscr{C}(m,k;\sigma
,-n+i+(j-1)\sigma)}{\sigma^{k}}.
\end{eqnarray*}
Moreover, given incomplete information on $(X_{1},\ldots,X_{n})$, the
Bayesian nonparametric estimator coincides with
\begin{eqnarray*}
\mathcal{R}_{m}^{(n,j)}&=& j-\frac{1}{\mathscr{C}(n,j;\sigma)}\sum
_{s=1}^{n-(j-1)}{n\choose s}\mathscr{C}(s,1;\sigma)
\mathscr {C}(n-s,j-1;\sigma)
\\
&&\hphantom{ j-}{}
\times\sum_{k=0}^{m}
\frac{V_{n+m,j+k}}{V_{n,j}}\frac
{\mathscr{C}(m,k;\sigma,-n+s+(j-1)\sigma)}{\sigma^{k}}.
\end{eqnarray*}
Here $m_{i}\geq0$ denotes the number of distinct species observed in
the initial sample with frequency~$i$.
\end{cor}

The distributions of $R_{m}^{(n,j,\mathbf{n})}$ and $R_{m}^{(n,j)}$,
under the Ewens--Pitman sampling model, are specified in the next
propositions. We devote special attention to the Ewens--Pitman sampling
model because it has proven suitable for inference in species sampling
problems, particularly in genomics. See, for example, \cite{Lij207}
and \cite{Fav12a} for details. The corresponding results for the
Ewens sampling model can be recovered by letting $\sigma\rightarrow0$
and applying equation 2.63 in \cite{Cha05}.
%
\begin{prp}\label{prop1_1}
Under the Ewens--Pitman sampling model, the distribution of
$R_{m}^{(n,j,\mathbf{n})}$ coincides with
%
\begin{eqnarray}
\label{eq:prob1_twopar} &&\P \bigl[R_{m}^{(n,j,\mathbf{n})}=x \bigr]
\nonumber\\
&&\quad
=\frac{1}{(\theta+n)_{m\uparrow1}}(-1)^{j}\sum
_{v=j-x}^{j}{v\choose j-x } (-1)^{v+x}
\\
&&\qquad {}
\times\sum_{\{c_{1},\ldots,c_{v}\}\in\mathcal
{C}_{j,v}}\Biggl(\theta+n-\sum
_{i=1}^{v}n_{c_{i}}+\sigma v
\Biggr)_{m\uparrow1}\nonumber\vspace*{-1pt}
\end{eqnarray}
and\vspace*{-1pt}
%
\begin{eqnarray}
\label{eq:est1_twopar} \mathcal{R}_{m}^{(n,j,\mathbf{n})}=j-\frac{1}{(\theta+n)_{m\uparrow
1}}\sum
_{i=1}^{n}m_{i}(\theta+n-i+
\sigma)_{m\uparrow1}.
\end{eqnarray}
The random variable $R_{m}^{(n,j, \mathbf{n})}$ assigns positive
probability to any integer value $x$ such that $0\leq x\le\min(j,m)$.
\end{prp}
%
\begin{prp}\label{prop1_2}
Under the Ewens--Pitman sampling model, the distribution of
$R_{m}^{(n,j)}$ coincides with
%
\begin{eqnarray}
\label{eq:prob2_twopar} &&\P \bigl[R_{m}^{(n,j)}=x \bigr]
\nonumber\\
&&\quad
=\frac{1}{\mathscr{C}(n,j;\sigma)(\theta+n)_{m\uparrow
1}}(-1)^{j}\sum
_{v=j-x}^{j}{v\choose j-x }(-1)^{v+x}
\\
&&\qquad {}
\times\sum_{s=v}^{n-(j-v)}{n\choose
s}(\theta+n-s+v\sigma )_{m\uparrow1}\mathscr{C}(s,v;\sigma)\nonumber
\mathscr{C}(n-s,j-v;\sigma)
\end{eqnarray}
and
%
\begin{eqnarray}
\label{eq:est2_twopar} &&\mathcal{R}_{m}^{(n,j)}=j-\frac{1}{\mathscr{C}(n,j;\sigma)(\theta
+n)_{m\uparrow1}}
\nonumber\\[-8pt]\\[-8pt]
&&\hphantom{j-}{}
\times\sum_{s=1}^{n-(j-1)}{n\choose
s }(\theta+n-s+\sigma )_{m\uparrow1}\mathscr{C}(s,1;\sigma)
\mathscr{C}(n-s,j-1;\sigma).\nonumber
\end{eqnarray}
The random variable $R_{m}^{(n,j)}$ assigns positive probability to any
integer value $x$ such that $0\leq x\le\min(j,m)$.\vadjust{\goodbreak}
\end{prp}

\subsection{Probabilities of re-observing old species with a certain
frequency}\label{subsec2}

In this section, we consider the distribution of the number of old
distinct species that are re-observed in $(X_{n+1},\ldots,X_{n+m})$
with frequency $0\leq l\leq m$, conditional on complete and incomplete
information on the random partition induced by the initial observed
sample $(X_{1},\ldots,X_{n})$. Note that the case $l=0$ is of
particular interest, representing the number of old distinct species
that are not re-observed in the additional sample. Formally, in the
context of complete information, we are interested in the random
variable $R_{l,m}^{(n,j,\mathbf{n})}$ which is defined in distribution as
%
\begin{equation}
\label{eq:def_main_2} \P\bigl[R_{l,m}^{(n,j,\mathbf{n})}=x\bigr]=\P \Biggl[\sum
_{i=1}^{K_{n}}\mathbh{1}_{\{S_{i}=l\}}=x \Bigl|
K_{n}=j,\mathbf {N}_{n}=\mathbf{n} \Biggr].
\end{equation}
In the context of incomplete information, we are interested in the
random variable $R_{l,m}^{(n,j)}$ which is defined in distribution as
%
\begin{equation}
\label{eq:def_main_2_part} \P\bigl[R_{l,m}^{(n,j)}=x\bigr]=\P \Biggl[\sum
_{i=1}^{K_{n}}\mathbh{1}_{\{
S_{i}=l\}}=x \Bigl|
K_{n}=j \Biggr].
\end{equation}
In the next theorem, we derive the factorial moments of the random
variables in \eqref{eq:def_main_2} and \eqref{eq:def_main_2_part}.
The factorial moment \eqref{teo2_eq1} is obtained by a direct
application of Theorem~1 in \cite{Fav12a}. With regards to the
factorial moment \eqref{teo2_eq2}, this is obtained from \eqref
{teo2_eq1} by suitably marginalizing the frequencies $\mathbf{N}_{n}$.
Again, these factorial moments lead to the corresponding distributions
by means of standard arguments involving probability generating functions.
%
\begin{thm}\label{teo2}
Let $(X_{i})_{i\geq1}$ be an exchangeable sequence directed by $\tilde
{P}_{G}$. Then, for any $0\leq l\leq m$ and any integer $r\geq1$ one has
%
\begin{eqnarray}
\label{teo2_eq1} &&\E \bigl[ \bigl(R_{l,m}^{(n,j,\mathbf{n})}
\bigr)_{r\downarrow
1} \bigr]
\nonumber\\
&&\quad
=r!{m\choose l,\ldots,l,m-rl}\sum_{\{c_{1},\ldots,c_{r}\}\in
\mathcal{C}_{j,r}}
\prod_{i=1}^{r}(n_{c_{i}}-
\sigma)_{l\uparrow1}
\\
&&\qquad {}
\times\sum_{k=0}^{m}
\frac{V_{n+m,j+k}}{V_{n,j}}\frac
{\mathscr{C}(m-rl,k;\sigma,-n+\sum_{i=1}^{r}n_{c_{i}}+(j-r)\sigma
)}{\sigma^{k}}\nonumber\vadjust{\goodbreak}
\end{eqnarray}
and
%
\begin{eqnarray}
\label{teo2_eq2} &&\E \bigl[ \bigl(R_{l,m}^{(n,j)}
\bigr)_{r\downarrow1} \bigr]
\nonumber\\
&&\quad
=\frac{r!}{\mathscr{C}(n,j;\sigma)}{m\choose l,\ldots ,l,m-rl}\bigl(-\sigma(1-
\sigma)_{(l-1)\uparrow1}\bigr)^{r}\\
&&\qquad {}
\times\sum_{s=r}^{n-(j-r)}{n\choose
s}\mathscr {C}(s,r;\sigma-l)\mathscr{C}(n-s,j-r;\sigma)
\nonumber\\
&&\qquad {}
\times\sum_{k=0}^{m}
\frac{V_{n+m,j+k}}{V_{n,j}}\frac
{\mathscr{C}(m-rl,k;\sigma,-n+s+(j-r)\sigma)}{\sigma^{k}},\nonumber
\end{eqnarray}
where $\mathcal{C}_{j,r}$ denotes the set of the $r$-combinations
(without repetitions) of the elements $\{1,\ldots,j\}$.
\end{thm}

Again, the distributions of $R_{l,m}^{(n,j,\mathbf{n})}$ and
$R_{l,m}^{(n,j)}$ are interpretable as the posterior distributions of
the number of old distinct species that are re-observed in
$(X_{n+1},\ldots,X_{n+m})$ with frequency $0\leq l\leq m$ given,
respectively, complete and incomplete information on the random
partition induced by $(X_{1},\ldots,X_{n})$. The corresponding
Bayesian nonparametric estimators, denoted by $\mathcal
{R}_{l,m}^{(n,j,\mathbf{n})}=\E[R_{l,m}^{(n,j,\mathbf{n})}]$ and
$\mathcal{R}_{l,m}^{(n,j)}=\E[R_{l,m}^{(n,j)}]$, are specified in the
next corollary. See Proposition~\ref{prop2_1} and Proposition~\ref
{prop2_2} for an expression for these estimators under the Ewens--Pitman
sampling model.
%
\begin{cor}
The Bayesian nonparametric estimator of the number of old distinct
species that are re-observed, with frequency $0\leq l\leq m$, in an
additional sample of size $m$, given complete information on
$(X_{1},\ldots,X_{n})$, coincides with
\begin{eqnarray*}
\mathcal{R}_{l,m}^{(n,j,\mathbf{n})}&=&{m\choose l}\sum
_{i=1}^{n}m_{i}(i-\sigma)_{l\uparrow1}
\\
&&{}
\times\sum_{k=0}^{m}
\frac{V_{n+m,j+k}}{V_{n,j}}\frac
{\mathscr{C}(m-l,k;\sigma,-n+i+(j-1)\sigma)}{\sigma^{k}}.
\end{eqnarray*}
Moreover, given incomplete information on $(X_{1},\ldots,X_{n})$ the
Bayesian nonparametric estimator coincides with
\begin{eqnarray*}
\mathcal{R}_{l,m}^{(n,j)}&=&\frac{1}{\mathscr{C}(n,j;\sigma
)}{m\choose l}
\bigl(-\sigma(1-\sigma)_{(l-1)\uparrow1}\bigr)
\\
&&{}
\times\sum_{s=1}^{n-(j-1)}{n\choose
s}\mathscr{C}(s,1;\sigma -l)\mathscr{C}(n-s,j-1;\sigma)
\\
&&{}
\times\sum_{k=0}^{m}
\frac{V_{n+m,j+k}}{V_{n,j}}\frac
{\mathscr{C}(m-l,k;\sigma,-n+s+(j-1)\sigma) }{\sigma^{k}}.
\end{eqnarray*}
Here $m_{i}\geq0$ denotes the number of distinct species observed in
the initial sample with frequency~$i$.
\end{cor}

Finally, the distributions of $R_{m}^{(n,j,\mathbf{n})}$ and
$R_{m}^{(n,j)}$, under the Ewens--Pitman sampling model, are specified
in the next propositions.
%
\begin{prp}\label{prop2_1}
Under the Ewens--Pitman sampling model, for any $0\leq l\leq m$, the
distribution of $R_{l,m}^{(n,j,\mathbf{n})}$ coincides with
%
\begin{eqnarray}
\label{eq:prob1_freq_twopar} &&\P \bigl[R_{l,m}^{(n,j,\mathbf{n})}=x \bigr]\nonumber
\\
&&\quad
=\frac{1}{(\theta+n)_{m\uparrow1}}\sum_{y=x}^{m}{y
\choose y-x}(-1)^{y-x}
\\
&&\qquad {}
\times{m\choose l,\ldots,l,m-yl}\sum_{\{c_{1},\ldots
,c_{y}\}\in\mathcal{C}_{j,y}}
\prod_{i=1}^{y}(n_{c_{i}}-\sigma
)_{l\uparrow1}\Biggl(\theta+n-\sum_{i=1}^{y}n_{c_{i}}+
\sigma y\Biggr)_{(m-yl)\uparrow1}\nonumber
\end{eqnarray}
and
%
\begin{eqnarray}
\label{eq:est1_freq_twopar} \mathcal{R}_{l,m}^{(n,j,\mathbf{n})}=\frac{1}{(\theta+n)_{m\uparrow
1}}{m
\choose l}\sum_{i=1}^{n}m_{i}(i-
\sigma)_{l\uparrow1}(\theta +n-i+\sigma)_{(m-l)\uparrow1}.
\end{eqnarray}
The random variable $R_{l,m}^{(n,j,\mathbf{n})}$ assigns positive
probability to any integer value $x$ such that $0\leq x\le\min(j,m)$.
\end{prp}
%
\begin{prp}\label{prop2_2}
Under the Ewens--Pitman sampling model, for any $0\leq l\leq m$, the
distribution of $R_{l,m}^{(n,j)}$ coincides with
%
\begin{eqnarray}
\label{eq:prob2_freq_twopar} &&\P \bigl[R_{l,m}^{(n,j)}=x \bigr]
\nonumber\\
&&\qquad
=\frac{1}{\mathscr{C}(n,j;\sigma)(\theta+n)_{m\uparrow
1}}\sum_{y=x}^{m}{y
\choose y-x}(-1)^{y-x}
\nonumber\\[-8pt]\\[-8pt]
&&\qquad {}
\times{m\choose l,\ldots,l,m-yl}\bigl(-\sigma(1-\sigma
)_{(l-1)\uparrow1}\bigr)^{y}\nonumber
\\
&&\qquad {}
\times\sum_{s=y}^{n-(j-y)}{n\choose
s}(\theta+n-s+\sigma y)_{(m-yl)\uparrow1}\mathscr{C}(s,y;\sigma-l)\mathscr
{C}(n-s,j-y;\sigma)\nonumber
\end{eqnarray}
and
%
\begin{eqnarray}
\label{eq:est2_freq_twopar} \mathcal{R}_{l,m}^{(n,j)}&=&\frac{1}{\mathscr{C}(n,j;\sigma)(\theta
+n)_{m\uparrow1}}{m
\choose l}\bigl(-\sigma(1-\sigma)_{(l-1)\uparrow1}\bigr)\nonumber
\\[-8pt]\\[-8pt]
&&{}
\times\sum_{s=1}^{n-(j-1)}{n\choose
s}(\theta+n-s+\sigma )_{(m-l)\uparrow1}\mathscr{C}(s,1;\sigma-l)\mathscr
{C}(n-s,j-1;\sigma).\nonumber
\end{eqnarray}
The random variable $R_{l,m}^{(n,j)}$ assigns positive probability to
any integer value $x$ such that $0\leq x\le\min(j,m)$.
\end{prp}

\subsection{Conditioning on almost-complete information}\label{subsec3}

We now consider the distribution of the number of old distinct species
that are re-observed in the additional sample $(X_{n+1},\ldots
,X_{n+m})$, conditional on almost-complete information. This
looking-backward species sampling problem can be seen as a
generalization of the problems discussed above. For any integer $p\in\{
1,\ldots,K_{n}\}$ let $\tau=\{\tau_{1},\ldots,\tau_{p}\}$ be a
collection of integers such that $1\leq\tau_{1}<\cdots<\tau_{p}\leq
K_{n}$ and define the subset of $p$ frequencies $\mathbf{N}_{\tau
,n}=(N_{\tau_{1}},\ldots, N_{\tau_{p}})$. In the context of
almost-complete information, we are interested in the random variables
$R_{m}^{(n,j,\mathbf{n}_{\tau})}$ and $R_{l,m}^{(n,j,\mathbf
{n}_{\tau})}$ which are defined in distribution as
%
\begin{equation}
\label{eq_totold_almostcomplete} \P\bigl[R_{m}^{(n,j,\mathbf{n}_{\tau})}=x\bigr]=\P \Biggl[\sum
_{i=1}^{K_{n}}\mathbh{1}_{\{S_{i}>0\}}=x \Bigl|
K_{n}=j,\mathbf{N}_{\tau
,n}=\mathbf{n}_{\tau} \Biggr]
\end{equation}
and
%
\begin{equation}
\label{eq_freqold_almostcomplete} \P\bigl[R_{l,m}^{(n,j,\mathbf{n}_{\tau})}=x\bigr]=\P \Biggl[\sum
_{i=1}^{K_{n}}\mathbh{1}_{\{S_{i}=l\}}=x \Bigl|
K_{n}=j,\mathbf{N}_{\tau
,n}=\mathbf{n}_{\tau} \Biggr].
\end{equation}
The following lemma is fundamental in determining the factorial moments
of the random variables introduced in \eqref{eq_totold_almostcomplete}
and \eqref{eq_freqold_almostcomplete} and, accordingly, to derive the
corresponding distributions.
%
\begin{lem}\label{margin_lem}
Let $(X_{i})_{i\geq1}$ be an exchangeable sequence directed by a
Gibbs-type random probability measure $\tilde{P}_{G}$. For any integer
$p\in\{1,\ldots,K_{n}\}$, denote by $\nu=\{\nu_{1},\ldots,\nu
_{K_{n}-p}\}$ the complement set of $\tau$ with $1\leq\nu
_{1}<\cdots<\nu_{K_{n}-p}\leq K_{n}$ and define the subset of
frequencies $\mathbf{N}_{\nu,n}:=(N_{\nu_{1}},\ldots,N_{\nu
_{K_{n}-p}})$. Then
%
\begin{eqnarray}
\label{eq:marginal_eppf} &&\P[\mathbf{N}_{\nu,n}=\mathbf{n}_{\nu} |
K_{n}=j,\mathbf {N}_{\tau,n}=\mathbf{n}_{\tau}]\nonumber
\\[-8pt]\\[-8pt]
&&\quad
=\frac{\sigma^{j-p}}{\mathscr{C}(n-\sum_{i=1}^{p}n_{\tau
_{i}},j-p;\sigma)}\frac{1}{(j-p)!}{n-\sum
_{i=1}^{p}n_{\tau
_{i}}\choose
n_{\nu_{1}},\ldots,n_{\nu_{j-p}}}\prod_{i=1}^{j-p}(1-
\sigma)_{(n_{\nu_{i}}-1)\uparrow1}.\nonumber
\end{eqnarray}
The random variable $\mathbf{N}_{\nu,n}=\mathbf{n}_{\nu} |
(K_{n}=j,\mathbf{N}_{\tau,n}=\mathbf{n}_{\tau})$ assigns positive
probability to the set $\mathcal{D}_{n-\sum_{i=1}^{p}n_{\tau_{i}},j-p}$.
\end{lem}
The factorial moments of $R_{m}^{(n,j,\mathbf{n}_{\tau})}$ and
$R_{l,m}^{(n,j,\mathbf{n}_{\tau})}$ are derived by means of Lemma~\ref{margin_lem} and along lines similar to the proof of Theorem~\ref
{teo1} and Theorem~\ref{teo2}, respectively. In particular, with
regard to the factorial moments of the random variables in \eqref
{eq_totold_almostcomplete}, one has
%
\begin{eqnarray}
\label{eq_totold_almostcomplete_mom} &&\E \bigl[\bigl(R_{m}^{(n)}\bigr)_{r\downarrow1}
| K_{n}=j,\mathbf{N}_{\tau
,n}=\mathbf{n}_{\tau} \bigr]\nonumber
\\
&&\quad
=\frac{r!}{\mathscr{C}(n-\sum_{i=1}^{p}n_{\tau
_{i}},j-p;\sigma)}\sum_{v_{1}=0}^{r}
\sum_{v_{2}=0}^{r}(-1)^{v_{1}+v_{2}}{j-v_{1}-v_{2}
\choose r-v_{1}-v_{2}}\nonumber
\\
&&\qquad {}
\times\sum_{\{d_{1},\ldots,d_{v_{1}}\}\in\mathcal
{C}_{p,v_{1}}}\sum
_{s=v_{2}}^{n-\sum_{i=1}^{p}n_{\tau
_{i}}-(j-p-v_{2})}{n-\sum_{i=1}^{p}n_{\tau_{i}}
\choose s}\\
&&\qquad {}\times\mathscr {C}(s,v_{2};\sigma)\mathscr{C}\Biggl(n-\sum
_{i=1}^{p}n_{\tau
_{i}}-s,j-p-v_{2};
\sigma\Biggr)\nonumber
\\
&&\qquad {}
\times\sum_{k=0}^{m}
\frac{V_{n+m,j+k}}{V_{n,j}}\frac
{\mathscr{C}(m,k;\sigma,-n+\sum_{i=1}^{v_{1}}n_{\tau
_{d_{i}}}+(j-v_{1}-v_{2})\sigma+s)}{\sigma^{k}}.\nonumber
\end{eqnarray}
We point out that \eqref{eq_totold_almostcomplete_mom} is a
generalization of both the results stated in Theorem~\ref{teo1}.
Indeed, by setting $\tau=j$ in \eqref{eq_totold_almostcomplete_mom}
one obtains \eqref{teo1_eq2}, whereas by setting $p=j$ in \eqref
{eq_totold_almostcomplete_mom} one obtains \eqref{teo1_eq1}. With
regard to the factorial moments of the random variables in \eqref
{eq_freqold_almostcomplete}, one has
%
\begin{eqnarray}
\label{eq_freqold_almostcomplete_mom} &&\E \bigl[\bigl(R_{l,m}^{(n)}
\bigr)_{r\downarrow1} | K_{n}=j,\mathbf{N}_{\tau
,n}=
\mathbf{n}_{\tau} \bigr]\nonumber
\\
&&\quad
=\frac{r!}{\mathscr{C}(n-\sum_{i=1}^{p}n_{\tau
_{i}},j-p;\sigma)}{m\choose l,\ldots,l,m-rl}\sum
_{v=0}^{r}\bigl(-\sigma (1-\sigma)_{(l-1)\uparrow1}
\bigr)^{r-v}\nonumber
\\
&&\qquad {}
\times\sum_{\{d_{1},\ldots,d_{v}\}\in\mathcal
{C}_{p,v}}\prod
_{i=1}^{v}(n_{\tau_{d_{i}}}-\sigma)_{l\uparrow1}\nonumber
\\[-8pt]\\[-8pt]
&&\qquad {}
\times\sum_{s=r-v}^{n-\sum_{i=1}^{p}n_{\tau
_{i}}-(j-p-(r-v))}{n-\sum
_{i=1}^{p}n_{\tau_{i}}\choose s}\nonumber\\
&&\qquad {}\times
\mathscr {C}(s,r-v;\sigma-l)\mathscr{C}\Biggl(n-\sum_{i=1}^{p}n_{\tau
_{i}}-s,j-p-(r-v);
\sigma\Biggr)\nonumber
\\
&&\qquad {}
\times\sum_{k=0}^{m}
\frac{V_{n+m,j+k}}{V_{n,j}}\frac
{\mathscr{C}(m-rl,k;\sigma,-n+\sum_{i=1}^{v}n_{\tau
_{d_{i}}}+s+(j-r)\sigma)}{\sigma^{k}}.\nonumber
\end{eqnarray}
Note that \eqref{eq_freqold_almostcomplete_mom} includes as special
cases both the results stated in Theorem~\ref{teo2}. Indeed, by
setting $\tau=j$ in \eqref{eq_freqold_almostcomplete_mom} one obtains
\eqref{teo2_eq2}, whereas by setting $p=j$ in \eqref
{eq_freqold_almostcomplete_mom} one obtains \eqref{teo2_eq1}.

\section{Numerical illustrations}\label{sec4}

We can now apply the derived conditional results which are
interpretable, from a Bayesian nonparametric standpoint, as estimators
or predictions. The range of problems to be addressed can be delineated
using the following hypothetical setting.
A nineteenth century naturalist samples a number of marine species in
an expedition to a remote island, reporting in his notebook the number
of distinct species sampled and their frequencies. We are interested in
estimating the abundance of a particular species observed at that point
in time. If all the data in the notebook are available, the
looking-backward estimators of Theorems \ref{teo1} and \ref{teo2}
which condition on complete information can be applied to solve this
problem. Now suppose that certain critical pages of the notebook are
missing, and the only datum available is the number of distinct species
in a sample of known size. This corresponds to the setting of
incomplete information.

In a general application, the species could be words in a text,
mutations of a gene in a population, or the names of newborns in a year.
The availability of complete or incomplete information could be
determined by constraints of the experimental method used or, in the
case of a meta-analysis, restrictions of access to data. For example,
techniques routinely used in biology provide indications about presence
or absence of a particular species, say a particular
bacterium or a genetic mutation of interest, but are not suitable for
measuring the relative species abundance. The experimental techniques,
in these cases,
produce datasets with partial information.

\begin{figure}

\includegraphics{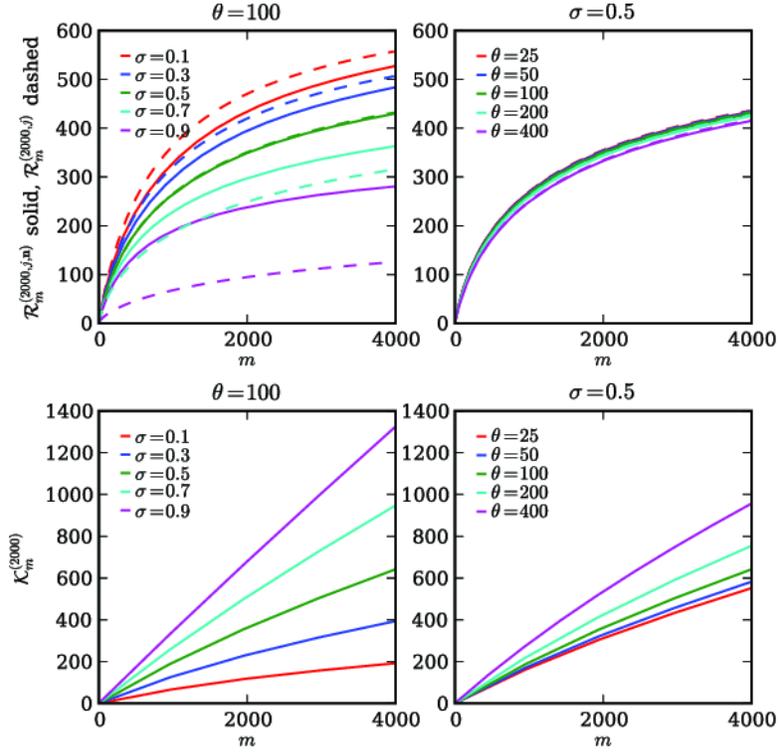}

\caption{Estimators for the number of old and new distinct species
observed as a function of the size $m$ of the additional sample. An
initial sample of $n=2000$ steps was drawn from the Ewens--Pitman
sampling model with $\theta=100$ and $\sigma=0.5$. The top panels
show estimators for the number of old species under complete
information, $\mathcal{R}^{(2000,j,\mathbf{n})}_m$, and incomplete
information, $\tilde{\mathcal{R}}^{(2000,j)}_m$. The bottom panels
show the estimator $\mathcal{K}^{(2000)}_m$ for the number of new
species. The panels on the left show estimators computed under $\theta
=100$ and allowing $\sigma$ to vary. The panels on the right show
estimators computed under $\sigma=0.5$ and allowing $\theta$ to vary.}
\label{fig:prima}
\end{figure}

We illustrate an application of the derived looking-backward estimators
in a simulation study. Two thousand samples were simulated from the
Ewens--Pitman sampling model with $\theta=100$ and $\sigma=0.5$. The
top row of panels in Figure~\ref{fig:prima} show the conditional
expectations of the number of re-observed species in an additional
experiment with sample sizes ranging from 0 to\vadjust{\goodbreak} 4000. These two panels
display discrepancies of the estimates under complete versus partial
information and illustrate sensitivity to the choice of the parameters
$\theta$ and $\sigma$. The estimates were computed across a range of
possible prior parameters, including the true data distribution.
Interestingly, the divergence between the two estimators depends more
heavily on $\sigma$ and is minimized when the parameter match those of
the true data distribution. We refer to \cite{Lij207} for detailed
arguments on practical selection of the prior parameters in this model.
The second row of panels, in contrast, displays estimates for the
number of new species in the additional sample. In this case the
estimates are identical under complete and partial information.
%

\begin{figure}

\includegraphics{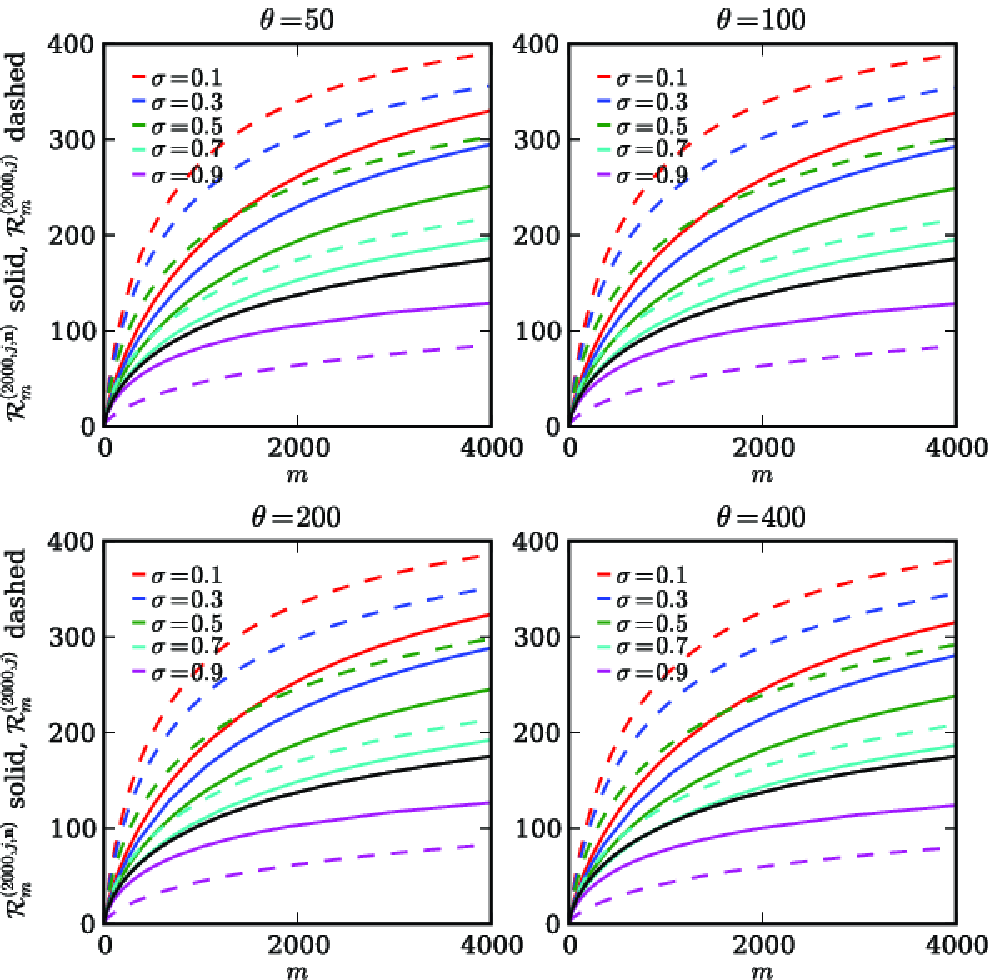}

\caption{Estimators for the number of old distinct species observed as
a function of the size $m$ of the additional sample. An initial sample
of $n=2000$ steps was drawn from a zeta distribution with scale
parameter $1.3$. Each panel shows estimators computed with a fixed
$\theta$ and allowing $\sigma$ to vary. The black line in each figure
shows the expected number of old distinct species in the sampling model.}
\label{fig:seconda}
\end{figure}

Figure~\ref{fig:seconda} considers simulated data that have not been
sampled from the Ewens--Pitman sampling model. Here, the sample was
generated from a Zeta distribution, whose power law behavior is common
in applications, and analyses were still performed using the
Ewens--Pitman sampling model. Looking-backward estimators under complete
and incomplete information are displayed for several prior parameters values.
These are consistent with the relationship between the choice of the
model parameters and the resulting conditional expectations shown in
Figure~\ref{fig:prima}.
Figure~\ref{fig:seconda} also displays (black line) the conditional
expectations under the true zeta sampling model, assumed unknown to the
investigator.
%

\begin{figure}

\includegraphics{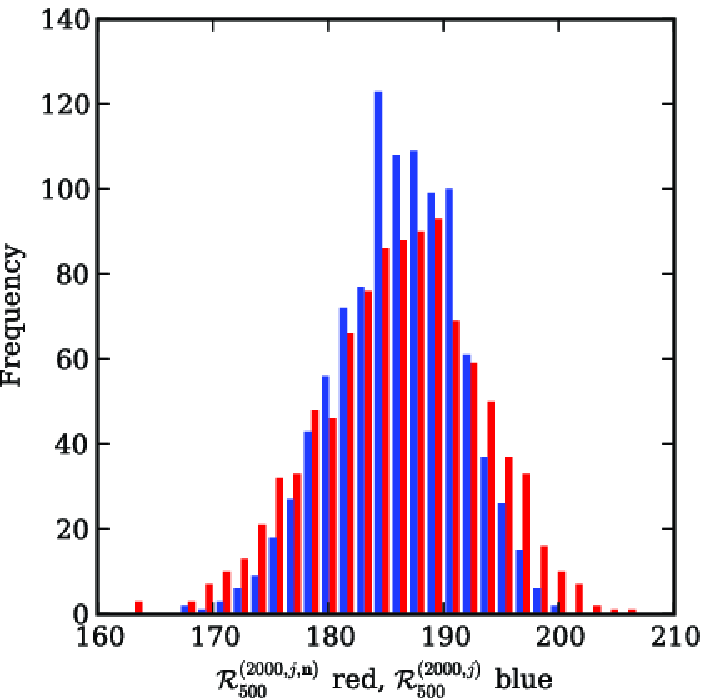}

\caption{Histograms of the estimators for the number of old species
under complete information, $\mathcal{R}^{(2000,j,\mathbf
{n})}_{500}$, and incomplete information, $\tilde{\mathcal
{R}}^{(2000,j)}_{500}$. To construct the histograms, these estimators
were computed conditional on 1000 independent initial samples of length
$n=2000$ each, which were drawn from the Ewens--Pitman sampling model
with $\theta=100$ and $\sigma=0.5$.}
\label{fig:terza}
\end{figure}

The simulations in Figure~\ref{fig:prima} were iterated, generating
1000 independent datasets of size $n=2000$ from the Ewens--Pitman
sampling model with $\theta=100$ and $\sigma=0.5$. Figure~\ref{fig:terza} shows the distribution of the estimator for the number of
distinct old species re-observed in an additional sample of size 500.
The blue and red histograms correspond to the estimator under complete
and incomplete information, respectively. As expected, the estimators
have the same mean but the estimator fit to complete information has
slightly higher variance.
%

\begin{appendix}

\section*{Appendix}\label{app}

\subsection{Proofs of the results in Section \texorpdfstring{\protect\ref{subsec1}}{3.1}}

\begin{pf*}{Proof of Theorem~\ref{teo1}} With regard to the $r$th
factorial moment of $R_{m}^{(n,j,\mathbf{n})}$, this is obtained by a
direct application of Theorem~1 in \cite{Fav12a}. Indeed, by means
of the Vandermonde's identity one has
%
\begin{eqnarray}
\label{eq:cond_mom_1}  \bigl(R^{(n,j,\mathbf{n})}_{m} \bigr)_{r\downarrow1}=\sum
_{v=0}^{r}{r\choose v}(-1)^{v}(j-v)_{(r-v)\downarrow1}
\bigl(R_{0,m}^{(n,j,\mathbf{n})} \bigr)_{v\downarrow1}.
\end{eqnarray}
Theorem~1 in \cite{Fav12a} then leads to \eqref{teo1_eq1} by
taking the expected value of both sides of \eqref{eq:cond_mom_1}. This
completes the first part of the proof. With regard to $r$th factorial
moment of the random variable $R_{m}^{(n,j)}$, by combining \eqref
{teo1_eq1} with the distributions displayed in \eqref{eq:gibbs_defi}
and \eqref{eq:gibbs_dist}, we write
%
\begin{eqnarray}
\label{eq:start} &&\E \bigl[ \bigl(R_{m}^{(n,j)}
\bigr)_{r\downarrow1} \bigr]\nonumber
\\
&&\quad
=\frac{\sigma^{j}}{\mathscr{C}(n,j;\sigma)} \sum_{v=0}^{r}{r
\choose v}(-1)^{v}(j-v)_{(r-v)\downarrow1}
\\
&&\qquad {}
\times\frac{1}{j!}\sum_{(n_{1},\ldots,n_{j})\in\mathcal
{D}_{n,j}}{n
\choose n_{1},\ldots,n_{j}}\prod
_{i=1}^{j}(1-\sigma )_{(n_{i}-1)\uparrow1}\nonumber
\\
&&\qquad {}
\times v!\sum_{\{c_{1},\ldots,c_{v}\}\in\mathcal
{C}_{j,v}}\sum
_{k=0}^{m}\frac{V_{n+m,j+k}}{V_{n,j}}\frac{\mathscr
{C}(m,k;\sigma,-n+\sum_{i=1}^{v}n_{c_{i}}+(j-v)\sigma)}{\sigma^{k}}\nonumber
\end{eqnarray}
and prove that it coincides with \eqref{teo1_eq2}. The proof is mainly
devoted to solve the sums over the indexes $n_{1},\ldots,n_{j}$ and
$c_{1},\ldots,c_{v}$. Once these sums are solved, then \eqref
{teo1_eq2} follows by some algebra involving factorial numbers and
noncentral generalized factorial coefficients. By means of equation
2.61 in \cite{Cha05}, and using the fact that $\mathcal{C}_{j,v}$
has cardinality ${j\choose v}$, from \eqref{eq:start} one has
%
\begin{eqnarray}
\label{eq:start_1} &&\E \bigl[ \bigl(R_{m}^{(n,j)}
\bigr)_{r\downarrow1} \bigr]\nonumber
\\
&&\quad
=\frac{\sigma^{j}}{\mathscr{C}(n,j;\sigma)}\sum_{k=0}^{m}
\frac{V_{n+m,j+k}}{V_{n,j}}\frac{1}{\sigma^{k}}\sum_{v=0}^{r}{r
\choose v}(-1)^{v}(j-v)_{(r-v)\downarrow1}\nonumber
\\
&&\qquad {}
\times\sum_{s_{1}=1}^{n-j+1}\sum
_{s_{2}=1}^{n-j+1-(s_{1}-1)}\cdots\sum
_{s_{v}=1}^{n-j+1-\sum
_{i=1}^{v-1}(s_{i}-1)}{n\choose s_{1},
\ldots,s_{v},n-\sum_{i=1}^{v}s_{i}}\\
&&\qquad {}\times
\prod_{i=1}^{v}(1-\sigma)_{(s_{i}-1)\uparrow1}\nonumber
\\
&&\qquad {}
\times\frac{1}{\sigma^{j-v}}\mathscr{C}\Biggl(m,k;\sigma ,-n+\sum
_{i=1}^{v}s_{i}+(j-v)\sigma
\Biggr)\mathscr{C}\Biggl(n-\sum_{i=1}^{v}s_{i},j-v;
\sigma\Biggr).\nonumber
\end{eqnarray}
In order to solve the nested sums over the indexes $s_{1},\ldots
,s_{v}$ in \eqref{eq:start_1}, we first deal with the sum over the
index $s_{v}$ and then we introduce a suitable recursive argument for
solving the remaining sums over the indexes $s_{1},\ldots,s_{v-1}$.
First, recall that for any $x\geq0$ and $0\leq y\leq x $, for any
$a>0$, $b>0$, $c>0$ and for any real number $d$ one has the following identity
%
\begin{equation}
\label{eq:gen_fact_coeff} {y+c\choose y}\mathscr{C}(x,y+c;d,a+b)=\sum
_{j=y}^{x-c}{x\choose j}\mathscr{C}(j,y;d,a)
\mathscr{C}(x-j,c;d,b).
\end{equation}
See Chapter~2 of \cite{Cha05} for details. Then, let us consider
the sum over the index $s_{v}$ in \eqref{eq:start_1}, that is,
\begin{eqnarray*}
&&\sum_{s_{v}=1}^{n-j+1-\sum_{i=1}^{v-1}(s_{i}-1)}{n\choose
s_{1},\ldots,s_{v},n-\sum_{i=1}^{v}s_{i}}
\prod_{i=1}^{v}(1-\sigma )_{(s_{i}-1)\uparrow1}
\\
&&\qquad {} \times\frac{1}{\sigma^{j-v}}\mathscr{C}\Biggl(m,k;\sigma,-n+\sum
_{i=1}^{v}s_{i}+(j-v)\sigma\Biggr)
\mathscr{C}\Biggl(n-\sum_{i=1}^{v}s_{i},j-v;
\sigma\Biggr)
\\
&&\quad  =\frac{1}{\sigma^{j-v}}{n\choose s_{1},\ldots,s_{v-1},n-\sum
_{i=1}^{v-1}s_{i}}\prod
_{i=1}^{v-1}(1-\sigma)_{(s_{i}-1)\uparrow1}
\\
&&\qquad {} \times\sum_{s_{v}=1}^{n-j+1-\sum_{i=1}^{v-1}(s_{i}-1)}{n-\sum
_{i=1}^{v-1}s_{i}\choose
s_{v}}(1-\sigma)_{(s_{v}-1)\uparrow1}
\\
&&\qquad {} \times\mathscr{C}\Biggl(m,k;\sigma,-n+\sum_{i=1}^{v-1}s_{i}+s_{v}+(j-v)
\sigma\Biggr)\mathscr{C}\Biggl(n-\sum_{i=1}^{v-1}s_{i}-s_{v},j-v;
\sigma\Biggr).
\end{eqnarray*}
By a direct application of \eqref{eq:gen_fact_coeff} to the
coefficients $\mathscr{C}(m,k;\sigma,-n+\sum_{i=1}^{v}s_{i}+(j-v)\sigma)$ and $\mathscr{C}(n-\sum_{i=1}^{v}s_{i},j-v;\sigma)$ we can write the last expression in the
following expanded form
\begin{eqnarray*}
&&\frac{1}{\sigma^{j-v}}{n\choose s_{1},\ldots,s_{v-1},n-\sum
_{i=1}^{v-1}s_{i}}\prod
_{i=1}^{v-1}(1-\sigma)_{(s_{i}-1)\uparrow1}
\\[3pt]
&&\qquad {} \times\sum_{t=k}^{m}{m\choose t}
\mathscr{C}(t,k;\sigma)\sum_{l=j-v}^{n-\sum_{i=1}^{v-1}s_{i}-1}{n-
\sum_{i=1}^{v-1}s_{i}\choose l}
\mathscr{C}\bigl(l,j-v;\sigma,-(j-v)\sigma\bigr)
\\[3pt]
&&\qquad {} \times\sum_{s_{v}=1}^{n-l-\sum_{i=1}^{v-1}s_{i}}{n-l-\sum
_{i=1}^{v-1}s_{i}\choose
s_{v}}(1-\sigma)_{(s_{v}-1)\uparrow1}
\\[3pt]
&&\qquad {} \times\Biggl(n-\sum_{i=1}^{v-1}s_{i}-s_{v}-(j-v)
\sigma\Biggr)_{(m-t)\uparrow
1}\bigl(-(j-v)\sigma\bigr)_{(n-l-\sum_{i=1}^{v-1}s_{i}-s_{v})\uparrow1}\vspace*{1pt}
\end{eqnarray*}
(by the Vandermonde's identity to expand $(n-
\sum_{i=1}^{v-1}s_{i}-s_{v}-(j-v)
\sigma)_{(m-t)\uparrow1}$)\vspace*{1pt}
\begin{eqnarray*}
&&\quad  =\frac{1}{\sigma^{j-v}}{n\choose s_{1},\ldots,s_{v-1},n-\sum
_{i=1}^{v-1}s_{i}}\prod
_{i=1}^{v-1}(1-\sigma)_{(s_{i}-1)\uparrow1}
\\[3pt]
&&\qquad {} \times\sum_{t=k}^{m}{m\choose t}
\mathscr{C}(t,k;\sigma)\sum_{h=0}^{m-t}{m-t
\choose h}\bigl(-(j-v)\sigma\bigr)_{h\uparrow1}
\\[3pt]
&&\qquad {} \times\sum_{l=j-v}^{n-\sum_{i=1}^{v-1}s_{i}-1}{n-\sum
_{i=1}^{v-1}s_{i}\choose
l}(l)_{(m-t-h)\uparrow1}\mathscr {C}\bigl(l,j-v;\sigma,-(j-v)\sigma\bigr)
\\[3pt]
&&\qquad {} \times\sum_{s_{v}=1}^{n-l-\sum_{i=1}^{v-1}s_{i}}{n-l-\sum
_{i=1}^{v-1}s_{i}\choose
s_{v}}\\[3pt]
&&\qquad\qquad\qquad\qquad {} \times(1-\sigma)_{(s_{v}-1)\uparrow
1}\bigl(-(j-v)\sigma+h
\bigr)_{(n-l-\sum_{i=1}^{v-1}s_{i}-s_{v})\uparrow1 }\vspace*{1pt}
\end{eqnarray*}
(by Equation 2.56 in \cite{Cha05} to solve the sum over the index
$s_{v}$)\vspace*{1pt}
\begin{eqnarray*}
&&\quad  =\frac{1}{\sigma^{j-v}}{n\choose s_{1},\ldots,s_{v-1},n-\sum
_{i=1}^{v-1}s_{i}}\prod
_{i=1}^{v-1}(1-\sigma)_{(s_{i}-1)\uparrow1}
\\[3pt]
&&\qquad {} \times\sum_{t=k}^{m}{m\choose t}
\mathscr{C}(t,k;\sigma)\sum_{h=0}^{m-t}{m-t
\choose h}\bigl(-(j-v)\sigma\bigr)_{h\uparrow1}
\\[3pt]
&&\qquad {} \times\sum_{l=j-v}^{n-\sum_{i=1}^{v-1}s_{i}-1}{n-\sum
_{i=1}^{v-1}s_{i}\choose
l}(l)_{(m-t-h)\uparrow1}\mathscr {C}\bigl(l,j-v;\sigma,-(j-v)\sigma\bigr)
\\[3pt]
&&\qquad {} \times\frac{1}{\sigma}\mathscr{C}\Biggl(n-l-\sum
_{i=1}^{v-1}s_{i},1;\sigma,(j-v)\sigma-h
\Biggr)\vspace*{1pt}
\end{eqnarray*}
providing the solution for the innermost nested sum over the index
$s_{v}$. Therefore, according to the last identity,\vadjust{\goodbreak} the $r$th factorial
moment of $R_{m}^{(n,j)}$ has the following reduced expression
%
\begin{eqnarray}
\label{eq:start_2} &&\E \bigl[ \bigl(R_{m}^{(n,j)}
\bigr)_{r\downarrow1} \bigr]\nonumber
\\
&&\quad
=\frac{\sigma^{j}}{\mathscr{C}(n,j;\sigma)}\sum_{k=0}^{m}
\frac{V_{n+m,j+k}}{V_{n,j}}\frac{1}{\sigma^{k}}\sum_{v=0}^{r}{r
\choose v}(-1)^{v}(j-v)_{(r-v)\downarrow1}\nonumber
\\
&&\qquad {}
\times\sum_{s_{1}=1}^{n-j+1}\sum
_{s_{2}=1}^{n-j+1-(s_{1}-1)}\cdots\sum
_{s_{v-1}=1}^{n-j+1-\sum
_{i=1}^{v-2}(s_{i}-1)}{n\choose s_{1},
\ldots,s_{v-1},n-\sum_{i=1}^{v-1}s_{i}}\nonumber\\
&&\qquad {}\times\prod_{i=1}^{v-1}(1-\sigma)_{(s_{i}-1)\uparrow1}
\\
&&\qquad {}
\times\sum_{t=k}^{m}{m\choose
t}\mathscr{C}(t,k;\sigma )\sum_{h=0}^{m-t}{m-t
\choose h}\bigl(-(j-v)\sigma\bigr)_{h\uparrow1}\nonumber
\\
&&\qquad {}
\times\sum_{l=j-v}^{n-\sum_{i=1}^{v-1}s_{i}-1}{n-\sum
_{i=1}^{v-1}s_{i}\choose
l}(l)_{(m-t-h)\uparrow1}\mathscr {C}\bigl(l,j-v;\sigma,-(j-v)\sigma\bigr)\nonumber
\\
&&\qquad {}
\times\frac{1}{\sigma^{j-v+1}}\mathscr{C}\Biggl(n-l-\sum
_{i=1}^{v-1}s_{i},1;\sigma,(j-v)\sigma-h
\Biggr).\nonumber
\end{eqnarray}
Starting from \eqref{eq:start_2} we can now introduce a recursive
argument to solve the remaining nested sums over the indexes
$s_{1},\ldots,s_{v-1}$. In particular, consider the sum over the index
$s_{v-1}$, that~is,
%
\begin{eqnarray}
\label{eq:second_nested} &&\sum_{s_{v-1}=1}^{n-j+1-\sum_{i=1}^{v-2}(s_{i}-1)}{n\choose
s_{1},\ldots,s_{v-1},n-\sum_{i=1}^{v-1}s_{i}}
\prod_{i=1}^{v-1}(1-\sigma)_{(s_{i}-1)\uparrow1}\nonumber
\\
&&\quad {}
\times\sum_{t=k}^{m}{m\choose t}
\mathscr{C}(t,k;\sigma)\sum_{h=0}^{m-t}{m-t
\choose h}\bigl(-(j-v)\sigma\bigr)_{h\uparrow1}\nonumber
\\[-8pt]\\[-8pt]
&&\quad {}
\times\sum_{l=j-v}^{n-\sum_{i=1}^{v-1}s_{i}-1}{n-\sum
_{i=1}^{v-1}s_{i}\choose
l}(l)_{(m-t-h)\uparrow1}\mathscr {C}\bigl(l,j-v;\sigma,-(j-v)\sigma\bigr)\nonumber
\\
&&\quad {}
\times\frac{1}{\sigma^{j-v+1}}\mathscr{C}\Biggl(n-l-\sum
_{i=1}^{v-1}s_{i},1;\sigma,(j-v)\sigma-h
\Biggr)\nonumber
\end{eqnarray}
which can be written as
\begin{eqnarray*}
&&\frac{1}{\sigma^{j-v+1}}{n\choose s_{1},\ldots,s_{v-2},n-\sum
_{i=1}^{2}s_{i}}\prod
_{i=1}^{v-2}(1-\sigma)_{(s_{i}-1)\uparrow1}
\\
&&\qquad {} \times\sum_{t=k}^{m}{m\choose t}
\mathscr{C}(t,k;\sigma)\sum_{h=0}^{m-t}{m-t
\choose h}\bigl(-(j-v)\sigma\bigr)_{h\uparrow1}
\\
&&\qquad {} \times\sum_{l=j-v+1}^{n-\sum
_{i=1}^{v-2}s_{i}-1}(l-1)_{(m-t-h)\uparrow1}
\mathscr {C}\bigl(l-1,j-v;\sigma,-(j-v)\sigma\bigr)
\\
&&\qquad {} \times\sum_{s_{v-1}=1}^{n-l-\sum_{i=1}^{v-2}s_{i}}{n-\sum
_{i=1}^{v-2}s_{i}\choose s_{v-1}}
{n-\sum_{i=1}^{v-2}s_{i}-s_{v-1}
\choose l-1}
\\
&&\qquad {} \times(1-\sigma)_{(s_{v-1}-1)\uparrow1}\mathscr{C}\Biggl(n-l+1-\sum
_{i=1}^{v-2}s_{i}-s_{v-1},1;
\sigma,(j-v)\sigma-h\Biggr)
\end{eqnarray*}
(by \eqref{eq:gen_fact_coeff} to expand
$\mathscr {C}(n-l-\sum_{i=1}^{v-2}s_{i}-s_{v-1}+1,1;
\sigma,(j-v)\sigma-h)$)
\begin{eqnarray*}
&&\quad  =\frac{1}{\sigma^{j-v+1}}{n\choose s_{1},\ldots,s_{v-2},n-\sum
_{i=1}^{2}s_{i}}\prod
_{i=1}^{v-2}(1-\sigma)_{(s_{i}-1)\uparrow1}
\\
&&\qquad {} \times\sum_{t=k}^{m}{m\choose t}
\mathscr{C}(t,k;\sigma)\sum_{h=0}^{m-t}{m-t
\choose h}\bigl(-(j-v)\sigma\bigr)_{h\uparrow1}
\\
&&\qquad {} \times\sum_{l=j-v+1}^{n-\sum
_{i=1}^{v-2}s_{i}-1}(l-1)_{(m-t-h)\uparrow1}
\mathscr {C}\bigl(l-1,j-v;\sigma,-(j-v)\sigma\bigr)
\\[1.5pt]
&&\qquad {} \times\sum_{z=1}^{n-l-\sum_{i=1}^{v-2}s_{i}}{n-\sum
_{i=1}^{v-2}s_{i}\choose l-1,z,n-l+1-z-
\sum_{i=1}^{v-2}s_{i}}\mathscr
{C}(z,1;\sigma)
\\[1pt]
&&\qquad {} \times\sum_{s_{v-1}=1}^{n-l+1-z-\sum
_{i=1}^{v-2}s_{i}}{n-l+1-z-\sum
_{i=1}^{v-2}s_{i}\choose
s_{v-1}}
\\[1pt]
&&\qquad {} \times(1-\sigma)_{(s_{v-1}-1)\uparrow1}\bigl(-(j-v)\sigma +h\bigr)_{(n-l+1-z-\sum_{i=1}^{v-2}s_{i}-s_{v-1})\uparrow1}
\end{eqnarray*}
(by equation 2.56 in \cite{Cha05} to solve the sum over the index
$s_{v-1}$)\vspace*{1.5pt}
\begin{eqnarray*}
&&\quad  =\frac{1}{\sigma^{j-v+1}}{n\choose s_{1},\ldots,s_{v-2},n-\sum
_{i=1}^{2}s_{i}}\prod
_{i=1}^{v-2}(1-\sigma)_{(s_{i}-1)\uparrow1}
\\[1pt]
&&\qquad {} \times\sum_{t=k}^{m}{m\choose t}
\mathscr{C}(t,k;\sigma)\sum_{h=0}^{m-t}{m-t
\choose h}\bigl(-(j-v)\sigma\bigr)_{h\uparrow1}
\\[1.5pt]
&&\qquad {} \times\sum_{l=j-v+1}^{n-\sum
_{i=1}^{v-2}s_{i}-1}(l-1)_{(m-t-h)\uparrow1}
\mathscr {C}\bigl(l-1,j-v;\sigma,-(j-v)\sigma\bigr)
\\[1pt]
&&\qquad {} \times\sum_{z=1}^{n-l-\sum_{i=1}^{v-2}s_{i}}{n-\sum
_{i=1}^{v-2}s_{i}\choose l-1,z,n-l+1-z-
\sum_{i=1}^{v-2}s_{i}}\mathscr
{C}(z,1;\sigma)
\\[1pt]
&&\qquad {} \times\frac{1}{\sigma}\mathscr{C}\Biggl(n-l+1-z-\sum
_{i=1}^{v-2}s_{i},1;\sigma,(j-v)\sigma-h
\Biggr)
\end{eqnarray*}
(by \eqref{eq:gen_fact_coeff} to solve the sum
over the index $z$)\vspace*{1.5pt}
\begin{eqnarray*}
&&\quad  ={2\choose1} {n\choose s_{1},\ldots,s_{v-2},n-\sum
_{i=1}^{2}s_{i}}\prod
_{i=1}^{v-2}(1-\sigma)_{(s_{i}-1)\uparrow1}
\\[1pt]
&&\qquad {} \times\sum_{t=k}^{m}{m\choose t}
\mathscr{C}(t,k;\sigma)\sum_{h=0}^{m-t}{m-t
\choose h}\bigl(-(j-v)\sigma\bigr)_{h\uparrow1}
\\
&&\qquad {} \times\sum_{l=j-v}^{n-\sum_{i=1}^{v-2}s_{i}}{n-\sum
_{i=1}^{v-2}s_{i}\choose
l}(l)_{(m-t-h)\uparrow1}\mathscr {C}\bigl(l,j-v;\sigma,-(j-v)\sigma\bigr)
\\
&&\qquad {} \times\frac{1}{\sigma^{j-v+2}}\mathscr{C}\Biggl(n-l-\sum
_{i=1}^{v-2}s_{i},2;\sigma,(j-v)\sigma-h
\Biggr).
\end{eqnarray*}
Note that the resulting expression has the same structure of the
summand in \eqref{eq:second_nested}. This fact suggests the
possibility of repeating the above arguments to each of the remaining
nested sums over the indexes $s_{v-2},\ldots,s_{1}$, respectively. In
particular, after a repeated application of these arguments we can
write the $r$th factorial moment of $R_{m}^{(n,j)}$ as follows
%
\begin{eqnarray}
\label{eq:start_3} &&\E \bigl[ \bigl(R_{m}^{(n,j)}
\bigr)_{r\downarrow1} \bigr]\nonumber
\\
&&\quad
=\frac{\sigma^{j}}{\mathscr{C}(n,j;\sigma)}\sum_{k=0}^{m}
\frac{V_{n+m,j+k}}{V_{n,j}}\frac{1}{\sigma^{k}}r!\sum_{v=0}^{r}{j-v
\choose r-v}(-1)^{v}\nonumber
\\[-8pt]\\[-8pt]
&&\qquad {}
\times\sum_{t=k}^{m}{m\choose
t}\mathscr{C}(t,k;\sigma )\sum_{h=0}^{m-t}{m-t
\choose h}\bigl(-(j-v)\sigma\bigr)_{h\uparrow1}\nonumber
\\
&&\qquad {}
\times\frac{1}{\sigma^{j}}\sum_{l=j-v}^{n-v}{n
\choose l}(l)_{(m-t-h)\uparrow1}\mathscr{C}\bigl(l,j-v;\sigma,-(j-v)\sigma \bigr)
\mathscr{C}\bigl(n-l,v;\sigma,(j-v)\sigma-h\bigr).\nonumber
\end{eqnarray}
Finally, a direct application of \eqref{eq:gen_fact_coeff} to expand
$\mathscr{C}(n-l,v;\sigma,(j-v)\sigma-h)$ we can write \eqref
{eq:start_3} as
\begin{eqnarray*}
&&\E \bigl[ \bigl(R_{m}^{(n,j)} \bigr)_{r\downarrow1} \bigr]
\\
&&\quad  =\frac{\sigma^{j}}{\mathscr{C}(n,j;\sigma)}\sum_{k=0}^{m}
\frac
{V_{n+m,j+k}}{V_{n,j}}\frac{1}{\sigma^{k}} r!\sum_{v=0}^{r}{j-v
\choose r-v}(-1)^{v}\frac{1}{\sigma^{j}}\sum
_{s=v}^{n-(j-v)}{n\choose s}\mathscr{C}(s,v;\sigma)
\\
&&\qquad {} \times\sum_{t=k}^{m}{m\choose t}
\mathscr{C}(t,k;\sigma)\sum_{h=0}^{m-t}{m-t
\choose h}\bigl(-(j-v)\sigma\bigr)_{h\uparrow1}
\\
&&\qquad {} \times\sum_{l=j-v}^{n-s}{n-s\choose
s}(l)_{(m-t-h)\uparrow
1}\bigl(-(j-v)\sigma+h\bigr)_{(n-s-l)\uparrow1}\mathscr{C}
\bigl(l,j-v;\sigma ,-(j-v)\sigma\bigr)
\end{eqnarray*}
which leads to \eqref{teo1_eq2} by means of \eqref{eq:gen_fact_coeff}
and some standard algebra involving factorial numbers and noncentral
generalized factorial coefficients. This completes the second part of
the proof.
\end{pf*}
\begin{pf*}{Proof of Proposition~\ref{prop1_1}} By combining the $r$th
factorial moment of $R_{m}^{(n,j,\mathbf{n})}$ in Theorem~\ref{teo1}
with $V_{n,j}$ displayed in \eqref{eq:defi_ewpit} one has
%
\begin{eqnarray}
\label{eq:start_prop1} &&\E \bigl[ \bigl(R_{m}^{(n,j,\mathbf{n})}
\bigr)_{r\downarrow
1} \bigr]\nonumber
\\
&&\quad
=\frac{r!}{(\theta+n)_{m\uparrow1}}\sum_{v=0}^{r}{j-v
\choose r-v}(-1)^{v}\nonumber
\\
&&\qquad {}
\times\sum_{\{c_{1},\ldots,c_{v}\}\in\mathcal
{C}_{j,v}}\sum
_{k=0}^{m} \biggl(\frac{\theta}{\sigma}+j
\biggr)_{k\uparrow1}\mathscr{C}\Biggl(m,k;\sigma,-n+\sum
_{i=1}^{v}n_{c_{i}}+(j-v)\sigma\Biggr)
\\
&&\quad
=\frac{r!}{(\theta+n)_{m\uparrow1}}\sum_{v=0}^{r}{j-v
\choose r-v}(-1)^{v}\nonumber
\\
&&\qquad {}
\times\sum_{\{c_{1},\ldots,c_{v}\}\in\mathcal
{C}_{j,v}}\Biggl(\theta+n-\sum
_{i=1}^{v}n_{c_{i}}+\sigma v
\Biggr)_{m\uparrow1},\nonumber
\end{eqnarray}
where the last identity follows from equation 2.49 in \cite{Cha05}.
Accordingly, \eqref{eq:est1_twopar} follows from \eqref
{eq:start_prop1} by setting $r=1$. Regarding \eqref{eq:prob1_twopar},
an inversion of the generating function for the $r$th factorial moment
in \eqref{eq:start_prop1} leads to
%
\begin{eqnarray}
\label{eq:start_prop2} &&\P \bigl[R_{m}^{(n,j,\mathbf{n})}=x \bigr]
\nonumber\\
&&\quad
=\frac{1}{(\theta+n)_{m\uparrow1}}\sum_{y\geq0}
\frac
{1}{x!}\frac{\ddr^{x}}{\ddr t^{x} }(t-1)^{x+y}\biggl\vert
_{t=0}
\\
&
&\qquad {}\times\sum_{v=0}^{x+y}{j-v
\choose x+y-v}(-1)^{v}\sum_{\{
c_{1},\ldots,c_{v}\}\in\mathcal{C}_{j,v}}\Biggl(
\theta+n-\sum_{i=1}^{v}n_{c_{i}}+
\sigma v\Biggr)_{m\uparrow1},\nonumber
\end{eqnarray}
where
%
\begin{eqnarray*}
\frac{\ddr^{x}}{\ddr t^{x} }(t-1)^{x+y}\biggl\vert _{t=0}=(-1)^{y}(x+y)_{x\downarrow1}.
\end{eqnarray*}
The proof is then completed by means of standard algebra involving
factorial numbers and binomial coefficients. Specifically, since
${j-v\choose x+y-v}=0$ for any $y>j-x$ then \eqref{eq:start_prop2} can
be written as
\begin{eqnarray*}
&&\P \bigl[R_{m}^{(n,j,\mathbf{n})}=x \bigr]
\\
&&\quad  =\frac{1}{(\theta+n)_{m\uparrow1}}\sum_{y=0}^{j}(-1)^{y-x}{y
\choose y-x}\sum_{v=0}^{y}(-1)^{v}{j-v
\choose y-v}
\\
&&\qquad {} \times\sum_{\{c_{1},\ldots,c_{v}\}\in\mathcal{C}_{j,v}}\Biggl(\theta +n-\sum
_{i=1}^{v}n_{c_{i}}+\sigma v
\Biggr)_{m\uparrow1}
\\
&&\quad  =\frac{1}{(\theta+n)_{m\uparrow1}}(-1)^{x}\sum_{v=0}^{j}
\sum_{y=0}^{j-v}(-1)^{y}{j-v
\choose y} {y+v\choose x}
\\
&&\qquad {} \times\sum_{\{c_{1},\ldots,c_{v}\}\in\mathcal{C}_{j,v}}\Biggl(\theta +n-\sum
_{i=1}^{v}n_{c_{i}}+\sigma v
\Biggr)_{m\uparrow1}
\\
&&\quad  =\frac{1}{(\theta+n)_{m\uparrow1}}(-1)^{x}\sum_{v=0}^{j}
(-1)^{j-v}{v\choose x-j+v }
\\
&&\qquad {} \times\sum_{\{c_{1},\ldots,c_{v}\}\in\mathcal{C}_{j,v}}\Biggl(\theta +n-\sum
_{i=1}^{v}n_{c_{i}}+\sigma v
\Biggr)_{m\uparrow1}
\end{eqnarray*}
which leads to \eqref{eq:prob1_twopar} by means of standard algebraic
manipulations involving factorial numbers.
\end{pf*}
\begin{pf*}{Proof of Proposition~\ref{prop1_2}} A combination of the
$r$th factorial moment of $R_{m}^{(n,j)}$ in Theorem~\ref{teo1} with
$V_{n,j}$ displayed in \eqref{eq:defi_ewpit} leads to
%
\begin{eqnarray}
\label{eq:start_prop3} &&\E \bigl[ \bigl(R_{m}^{(n,j)}
\bigr)_{r\downarrow1} \bigr]\nonumber
\\
&&\quad
=\frac{r!}{\mathscr{C}(n,j;\sigma)(\theta+n)_{m\uparrow
1}}\sum_{v=0}^{r}{j-v
\choose r-v}(-1)^{v}\nonumber
\\
&&\qquad {}
\times\sum_{s=v}^{n-(j-v)}{n\choose
s}\mathscr {C}(s,v;\sigma)\mathscr{C}(n-s,j-v;\sigma)\nonumber
\\[-8pt]\\[-8pt]
&&\qquad {}
\times\sum_{k=0}^{m} \biggl(
\frac{\theta}{\sigma}+j \biggr)_{k\uparrow1}\mathscr{C}\bigl(m,k;\sigma,-n+s+(j-v)
\sigma\bigr)\nonumber
\\
&&\quad
=\frac{r!}{\mathscr{C}(n,j;\sigma)(\theta+n)_{m\uparrow
1}}\sum_{v=0}^{r}{j-v
\choose r-v}(-1)^{v}\nonumber
\\
&&\qquad {}
\times\sum_{s=v}^{n-(j-v)}{n\choose
s}(\theta+n-s+v\sigma )_{m\uparrow1}\mathscr{C}(s,v;\sigma)
\mathscr{C}(n-s,j-v;\sigma),\nonumber
\end{eqnarray}
where the last identity follows from equation 2.49 in \cite{Cha05}.
Accordingly, \eqref{eq:est2_twopar} follows from \eqref
{eq:start_prop3} by setting $r=1$. Regarding \eqref{eq:prob2_twopar},
an inversion of the generating function for the $r$th factorial moment
in \eqref{eq:start_prop3} leads to
%
\begin{eqnarray}
\label{eq:start_prop4} &&\P \bigl[R_{m}^{(n,j)}=x \bigr]\nonumber
\\
&&\quad
=\frac{1}{\mathscr{C}(n,j;\sigma)(\theta+n)_{m\uparrow
1}}\sum_{y\geq0}
\frac{1}{x!}\frac{\ddr^{x}}{\ddr t^{x}
}(t-1)^{x+y}\biggl\vert
_{t=0}\nonumber
\\[-8pt]\\[-8pt]
&&\qquad {}
\times\sum_{v=0}^{x+y}{j-v
\choose x+y-v}(-1)^{v}\nonumber
\\
&&\qquad {}
\times\sum_{s=v}^{n-(j-v)}{n\choose
s}(\theta+n-s+v\sigma )_{m\uparrow1}\mathscr{C}(s,v;\sigma)
\mathscr{C}(n-s,j-v;\sigma),\nonumber
\end{eqnarray}
where\vspace*{-0.5pt}
%
\begin{eqnarray*}
\frac{\ddr^{x}}{\ddr t^{x} }(t-1)^{x+y}\biggl\vert _{t=0}=(-1)^{y}(x+y)_{x\downarrow1}.
\end{eqnarray*}
The proof is then completed by means of standard algebra involving
factorial numbers and binomial coefficients. Specifically, since
${j-v\choose x+y-v}=0$ for any $y>j-x$ then \eqref{eq:start_prop4} can
be written as
\begin{eqnarray*}
&&\P \bigl[R_{m}^{(n,j)}=x \bigr]
\\[-0.5pt]
&&\quad  =\frac{1}{\mathscr{C}(n,j;\sigma)(\theta+n)_{m\uparrow1}}\sum_{y=0}^{j}(-1)^{y-x}{y
\choose y-x}\sum_{v=0}^{y}{j-v\choose
y-v}(-1)^{v}
\\[-0.5pt]
&&\qquad {} \times\sum_{s=v}^{n-(j-v)}{n\choose s}(
\theta+n-s+v\sigma )_{m\uparrow1}\mathscr{C}(s,v;\sigma)\mathscr{C}(n-s,j-v;
\sigma)
\\[-0.5pt]
&&\quad  =\frac{1}{\mathscr{C}(n,j;\sigma)(\theta+n)_{m\uparrow
1}}(-1)^{x}\sum_{v=0}^{j}
\sum_{y=0}^{j-v}(-1)^{y}{j-v
\choose y} {y+v\choose x}
\\[-0.5pt]
&&\qquad {} \times\sum_{s=v}^{n-(j-v)}{n\choose s}(
\theta+n-s+v\sigma )_{m\uparrow1}\mathscr{C}(s,v;\sigma)\mathscr{C}(n-s,j-v;
\sigma)
\\[-0.5pt]
&&\quad  =\frac{1}{\mathscr{C}(n,j;\sigma)(\theta+n)_{m\uparrow
1}}(-1)^{x}\sum_{v=0}^{j}(-1)^{j-v}{v
\choose x-j+v }
\\[-0.5pt]
&&\qquad {} \times\sum_{s=v}^{n-(j-v)}{n\choose s}(
\theta+n-s+v\sigma )_{m\uparrow1}\mathscr{C}(s,v;\sigma)\mathscr{C}(n-s,j-v;
\sigma)
\end{eqnarray*}
which leads to \eqref{eq:prob2_twopar} by means of standard algebraic
manipulations involving factorial numbers.
\end{pf*}
\subsection{Proofs of the results in Section \texorpdfstring{\protect\ref{subsec2}}{3.2}}

\begin{pf*}{Proof of Theorem~\ref{teo2}} With regard to the $r$th
factorial moment of $R_{l,m}^{(n,j,\mathbf{n})}$, this is obtained by
a direct application of Theorem~1 in \cite{Fav12a}. This completes
the first part of the proof. With regard the $r$th factorial moment of
$R_{l,m}^{(n,j)}$, this is obtained by combining \eqref{teo2_eq1} with
the distributions displayed in \eqref{eq:gibbs_defi} and \eqref
{eq:gibbs_dist}. Specifically, we can write the following expression
%
\begin{eqnarray}
\label{eq:starting1} &&\E \bigl[ \bigl(R_{l,m}^{(n,j)}
\bigr)_{r\downarrow1} \bigr]\nonumber
\\
&&\quad
=\frac{\sigma^{j}}{\mathscr{C}(n,j;\sigma)}r!{m\choose l,\ldots,l,m-rl}
\nonumber\\
&&\qquad {}
\times\frac{1}{j!}\sum_{(n_{1},\ldots,n_{j})\in\mathcal
{D}_{n,j}}{n
\choose n_{1},\ldots,n_{j}}\prod
_{i=1}^{j}(1-\sigma )_{(n_{i}-1)\uparrow1}
\\
&&\qquad {}
\times\sum_{\{c_{1},\ldots,c_{r}\}\in\mathcal
{C}_{j,r}}\prod
_{i=1}^{r}(n_{c_{i}}-\sigma)_{l\uparrow1}\nonumber
\\
&&\qquad {}
\times\sum_{k=0}^{m}
\frac{V_{n+m,j+k}}{V_{n,j}}\frac
{\mathscr{C}(m-rl,k;\sigma,-n+\sum_{i=1}^{r}n_{c_{i}}+(j-r)\sigma)
}{\sigma^{k}}\nonumber
\end{eqnarray}
and prove that it coincides with \eqref{teo2_eq2}. As in Theorem~\ref
{teo1} the main issue consists in solving the sums over the collection
of indexes $n_{1},\ldots,n_{j}$ and $c_{1},\ldots,c_{r}$. First, by
means of equation 2.61 in \cite{Cha05}\vadjust{\goodbreak} and using the fact that
$\mathcal{C}_{j,r}$ has cardinality ${j\choose r }$, from \eqref
{eq:starting1} one has
%
\begin{eqnarray}
\label{eq:starting2} &&\E \bigl[ \bigl(R_{l,m}^{(n,j)}
\bigr)_{r\downarrow1} \bigr]\nonumber
\\
&&\quad
=\frac{\sigma^{j}}{\mathscr{C}(n,j;\sigma)}{m\choose l,\ldots,l,m-rl}\sum
_{k=0}^{m}\frac{V_{n+m,j+k}}{V_{n,j}}\frac
{1}{\sigma^{k}}\nonumber
\\
&&\qquad {}
\times\sum_{s_{1}=1}^{n-j+1}\sum
_{s_{2}=1}^{n-j+1-(s_{1}-1)}\cdots\sum
_{s_{r}=1}^{n-j+1-\sum
_{i=1}^{r-1}(s_{i}-1)}{n\choose s_{1},
\ldots,s_{r},n-\sum_{i=1}^{r}s_{i}}
\\
&&\qquad {}
\times\prod_{i=1}^{r}(1-
\sigma)_{(s_{i}-1)\uparrow
1}(s_{i}-\sigma)_{l\uparrow1}\nonumber
\\
&&\qquad {}
\times\frac{1}{\sigma^{j-r}}\mathscr{C}\Biggl(m-rl,k;\sigma ,-n+\sum
_{i=1}^{r}s_{i}+(j-r)\sigma
\Biggr)\mathscr{C}\Biggl(n-\sum_{i=1}^{r}s_{i},j-r;
\sigma\Biggr).\nonumber
\end{eqnarray}
As in the proof of Theorem~\ref{teo1}, in order to solve the nested
sums over the indexes $s_{1},\ldots,s_{r}$ in \eqref{eq:start_1} we
first deal with the sum over the index $s_{r}$. Recall that for any
$x\geq0$ and $0\leq y\leq x $, for any $a>0$, $b>0$, $c>0$ and for any
real number $d$ one has the following identity
%
\begin{equation}
\label{eq:gen_fact_coeff1} {y+c\choose y}\mathscr{C}(x,y+c;d,a+b)=\sum
_{j=y}^{x-c}{x\choose j}\mathscr{C}(j,y;d,a)
\mathscr{C}(x-j,c;d,b).
\end{equation}
See Chapter~2 of \cite{Cha05} for details. Then, let us consider
the sum over the index $s_{r}$ in \eqref{eq:starting2}, that is,
\begin{eqnarray*}
&&\sum_{s_{r}=1}^{n-j+1-\sum_{r=1}^{r-1}(s_{i}-1)}{n\choose
s_{1},\ldots,s_{r},n-\sum_{i=1}^{r}s_{i}}
\prod_{i=1}^{r}(1-\sigma
)_{(s_{i}-1)\uparrow1}(s_{i}-\sigma)_{l\uparrow1}
\\
&&\qquad {} \times\frac{1}{\sigma^{j-r}}\mathscr{C}\Biggl(m-rl,k;\sigma,-n+\sum
_{i=1}^{r}s_{i}+(j-r)\sigma\Biggr)
\mathscr{C}\Biggl(n-\sum_{i=1}^{r}s_{i},j-r;
\sigma\Biggr)
\\
&&\quad  ={n\choose s_{1},\ldots,s_{r-1},n-\sum
_{i=1}^{r-1}s_{i}}\prod
_{i=1}^{r-1}(1-\sigma)_{(s_{i}-1)\uparrow1}(s_{i}-
\sigma)_{l\uparrow
1}
\\
&&\qquad {} \times\sum_{s_{r}=1}^{n-j+1-\sum_{i=1}^{r-1}(s_{i}-1)}{n-\sum
_{i=1}^{r-1}s_{i}\choose
s_{r}}(1-\sigma)_{(s_{r}-1)\uparrow
1}(s_{r}-
\sigma)_{l\uparrow1}
\\
&&\qquad {} \times\frac{1}{\sigma^{j-r}}\mathscr{C}\Biggl(m-rl,k;\sigma,-n+\sum
_{i=1}^{r-1}s_{i}+s_{r}+(j-r)
\sigma\Biggr)\mathscr{C}\Biggl(n-\sum_{i=1}^{r-1}s_{i}+s_{r},j-r;
\sigma\Biggr).
\end{eqnarray*}
By a direct application of \eqref{eq:gen_fact_coeff1} to the
coefficients $\mathscr{C}(m-rl,k;\sigma,-n+\sum_{i=1}^{r}s_{i}+(j-r)\sigma)$ and $\mathscr{C}(n-\sum_{i=1}^{r}s_{i},j-r;\sigma)$ we can write the last expression in the
following expanded form
\begin{eqnarray*}
&&\frac{1}{\sigma^{j-r}}{n\choose s_{1},\ldots,s_{r-1},n-\sum
_{i=1}^{r-1}s_{i}}\prod
_{i=1}^{r-1}(1-\sigma)_{(s_{i}-1)\uparrow
1}(s_{i}-
\sigma)_{l\uparrow1}
\\
&&\qquad {} \times\sum_{t=k}^{m-rl}{m-rl\choose t}
\mathscr{C}(t,k;\sigma)\sum_{z=j-r}^{n-\sum_{i=1}^{r-1}s_{i}-1}{n-
\sum_{i=1}^{r-1}s_{i}\choose z}
\mathscr{C}\bigl(z,j-r;\sigma,-(j-r)\sigma\bigr)
\\
&&\qquad {} \times\sum_{s_{r}=1}^{n-z-\sum_{i=1}^{r-1}s_{i}}{n-z-\sum
_{i=1}^{r-1}s_{i}\choose
s_{r}}(1-\sigma)_{(s_{r}-1)\uparrow
1}(s_{r}-
\sigma)_{l\uparrow1}
\\
&&\qquad {} \times\Biggl(n-\sum_{i=1}^{r-1}s_{i}-s_{r}-(j-r)
\sigma \Biggr)_{(m-rl-t)\uparrow1}\bigl(-(j-r)\sigma\bigr)_{(n-z-\sum
_{i=1}^{r-1}s_{i}-s_{r})\uparrow1}
\end{eqnarray*}
(by the Vandermonde's identity to expand $(n-
\sum_{i=1}^{r-1}s_{i}-s_{r}-(j-r)
\sigma)_{(m-rl-t)\uparrow1}$)\vspace*{1pt}
\begin{eqnarray*}
&&\quad  =\frac{1}{\sigma^{j-r}}{n\choose s_{1},\ldots,s_{r-1},n-\sum
_{i=1}^{r-1}s_{i}}\prod
_{i=1}^{r-1}(1-\sigma)_{(s_{i}-1)\uparrow
1}(s_{i}-
\sigma)_{l\uparrow1}
\\
&&\qquad {} \times\sum_{t=k}^{m-rl}{m-rl\choose t}
\mathscr{C}(t,k;\sigma)\sum_{h=0}^{m-rl-t}{m-rl-t
\choose h}\bigl(-(j-r)\sigma\bigr)_{h\uparrow1}
\\
&&\qquad {} \times\sum_{z=j-r}^{n-\sum_{i=1}^{r-1}s_{i}-1}{n-\sum
_{i=1}^{r-1}s_{i}\choose
z}(z)_{(m-rl-t-h)\uparrow1}\mathscr {C}\bigl(z,j-r;\sigma,-(j-r)\sigma\bigr)
\\
&&\qquad {} \times(1-\sigma)_{(l-1)\uparrow1}\sum_{s_{r}=1}^{n-z-\sum
_{i=1}^{r-1}s_{i}}{n-z-
\sum_{i=1}^{r-1}s_{i}\choose
s_{r}}
\\
&&\qquad {} \times(l-\sigma)_{s_{r}\uparrow1}\bigl(-(j-r)\sigma+h\bigr)_{(n-z-\sum
_{i=1}^{r-1}s_{i}-s_{r})\uparrow1}
\end{eqnarray*}
(by equation 2.60 in \cite{Cha05} to solve the sum over the index
$s_{r}$)
\begin{eqnarray*}
&&\quad  =\frac{1}{\sigma^{j-r}}{n\choose s_{1},\ldots,s_{r-1},n-\sum
_{i=1}^{r-1}s_{i}}\prod
_{i=1}^{r-1}(1-\sigma)_{(s_{i}-1)\uparrow
1}(s_{i}-
\sigma)_{l\uparrow1}
\\
&&\qquad {} \times\sum_{t=k}^{m-rl}{m-rl\choose t}
\mathscr{C}(t,k;\sigma)\sum_{h=0}^{m-rl-t}{m-rl-t
\choose h}\bigl(-(j-r)\sigma\bigr)_{h\uparrow1}
\\
&&\qquad {} \times\sum_{z=j-r}^{n-\sum_{i=1}^{r-1}s_{i}-1}{n-\sum
_{i=1}^{r-1}s_{i}\choose
z}(z)_{(m-rl-t-h)\uparrow1}\mathscr {C}\bigl(z,j-r;\sigma,-(j-r)\sigma\bigr)
\\
&&\qquad {} \times(1-\sigma)_{(l-1)\uparrow1}(-1)\mathscr{C}\Biggl(n-z-\sum
_{i=1}^{r-1}s_{i};1;\sigma-l,(j-r)\sigma-h
\Biggr)
\end{eqnarray*}
providing the solution for the innermost nested sum over the index
$s_{r}$. Therefore, according to the last identity, the $r$th factorial
moment of $R_{l,m}^{(n,j)}$ in \eqref{eq:starting2} has the following
reduced expression
%
\begin{eqnarray}
\label{eq:starting3} &&\E \bigl[ \bigl(R_{l,m}^{(n,j)}
\bigr)_{r\downarrow1} \bigr]\nonumber
\\
&&\quad
=\frac{\sigma^{j}}{\mathscr{C}(n,j;\sigma)}{m\choose l,\ldots,l,m-rl}\sum
_{k=0}^{m}\frac{V_{n+m,j+k}}{V_{n,j}}\frac
{1}{\sigma^{k}}\nonumber
\\
&&\qquad {}
\times\sum_{s_{1}=1}^{n-j+1}\sum
_{s_{2}=1}^{n-j+1-(s_{1}-1)}\cdots\sum
_{s_{r-1}=1}^{n-j+1-\sum
_{i=1}^{r-2}(s_{i}-1)}{n\choose s_{1},
\ldots,s_{r-1},n-\sum_{i=1}^{r-1}s_{i}}\nonumber
\\
&&\qquad {}
\times\prod_{i=1}^{r-1}(1-
\sigma)_{(s_{i}-1)\uparrow
1}(s_{i}-\sigma)_{l\uparrow1}
\\
&&\qquad {}
\times\sum_{t=k}^{m-rl}{m-rl
\choose t}\mathscr {C}(t,k;\sigma)\sum_{h=0}^{m-rl-t}{m-rl-t
\choose h}\bigl(-(j-r)\sigma \bigr)_{h\uparrow1}\nonumber
\\
&&\qquad {}
\times\sum_{z=j-r}^{n-\sum_{i=1}^{r-1}s_{i}-1}{n-\sum
_{i=1}^{r-1}s_{i}\choose
z}(z)_{m-rl-t-h}\mathscr{C}\bigl(z,j-r;\sigma ,-(j-r)\sigma\bigr)\nonumber
\\
&&\qquad {}
\times\frac{1}{\sigma^{j-r}}(1-\sigma)_{(l-1)\uparrow
1}(-1)\mathscr{C}
\Biggl(n-z-\sum_{i=1}^{r-1}s_{i};1;
\sigma-l,(j-r)\sigma-h\Biggr).\nonumber
\end{eqnarray}
Starting from \eqref{eq:starting3} we can now repeatedly apply
equation 2.60 in \cite{Cha05} to solve the remaining sums over the
indexes $s_{1},\ldots,s_{r-1}$, respectively, starting from the index
$s_{r-1}$ and proceeding backward to the index $s_{1}$. As an example,
consider the sum over the index $s_{r-1}$, that is,
%
\begin{eqnarray}
\label{eq_sum_r1} &&\sum_{s_{r-1}=1}^{n-j+1-\sum_{i=1}^{r-2}(s_{i}-1)}{n\choose
s_{1},\ldots,s_{r-1},n-\sum_{i=1}^{r-1}s_{i}}
\prod_{i=1}^{r-s}(1-\sigma)_{(s_{i}-1)\uparrow1}(s_{i}-
\sigma)_{l\uparrow
1}\nonumber
\\
&&\qquad {}
\times\sum_{t=k}^{m-rl}{m-rl
\choose t }\mathscr {C}(t,k;\sigma)\sum_{h=0}^{m-rl-t}{m-rl-t
\choose h}\bigl(-(j-r)\sigma \bigr)_{h\uparrow1}\nonumber
\\[-8pt]\\[-8pt]
&&\qquad {}
\times\sum_{z=j-r}^{n-\sum_{i=1}^{r-1}s_{i}-1}{n-\sum
_{i=1}^{r-1}s_{i}\choose
z}(z)_{(m-rl-t-h)\uparrow1}\mathscr {C}\bigl(z,j-r;\sigma,-(j-r)\sigma\bigr)
\nonumber\\
&&\qquad {}
\times\frac{1}{\sigma^{j-r}}(1-\sigma)_{(l-1)\uparrow
1}(-1)\mathscr{C}
\Biggl(n-z-\sum_{i=1}^{r-1}s_{i};1;
\sigma-l,(j-r)\sigma-h\Biggr)\nonumber
\end{eqnarray}
which can be written as
\begin{eqnarray*}
&&\frac{1}{\sigma^{j-r}}{n\choose s_{1},\ldots,s_{r-2},n-\sum
_{i=1}^{r-2}s_{i}}\prod
_{i=1}^{r-2}(1-\sigma)_{(s_{i}-1)\uparrow
1}(s_{i}-
\sigma)_{l\uparrow1}
\\
&&\qquad {} \times\sum_{t=k}^{m-rl}{m-rl\choose t }
\mathscr{C}(t,k;\sigma )\sum_{h=0}^{m-rl-t}{m-rl-t
\choose h}\bigl(-(j-r)\sigma\bigr)_{h\uparrow1}
\\
&&\qquad {} \times\sum_{z=j-r+1}^{n-\sum_{i=1}^{r-2}s_{i}-1}{n-\sum
_{i=1}^{r-2}s_{i}\choose
z-1}(z-1)_{(m-rl-t-h)\uparrow1}\mathscr {C}\bigl(z-1,j-r;\sigma,-(j-r)\sigma\bigr)
\\
& &\qquad {}\times\bigl((1-\sigma)_{(l-1)\uparrow1}\bigr)^{2}(-1)\sum
_{s_{r-1}=1}^{n-z-\sum_{i=1}^{r-2}s_{i}}{n-z+1-\sum_{i=1}^{r-2}s_{i}
\choose s_{r-1}}
\\
&&\qquad {} \times(l-\sigma)_{s_{r-1}\uparrow1}\mathscr{C}\Biggl(n-z+1-\sum
_{i=1}^{r-2}s_{i}-s_{r-1},1;
\sigma-l,(j-r)\sigma-h\Biggr)
\end{eqnarray*}
(by equation 2.60 in \cite{Cha05} to solve the sum over the index
$s_{r-1}$)
\begin{eqnarray*}
&&\quad  ={n\choose s_{1},\ldots,s_{r-2},n-\sum
_{i=1}^{r-2}s_{i}}\prod
_{i=1}^{r-2}(1-\sigma)_{(s_{i}-1)\uparrow1}(s_{i}-
\sigma)_{l\uparrow
1}
\\
&&\qquad {} \times\sum_{t=k}^{m-rl}{m-rl\choose t }
\mathscr{C}(t,k;\sigma )\sum_{h=0}^{m-rl-t}{m-rl-t
\choose h}\bigl(-(j-r)\sigma\bigr)_{h\uparrow1}
\\
&&\qquad {} \times\sum_{z=j-r+1}^{n-\sum_{i=1}^{r-2}s_{i}-1}{n-\sum
_{i=1}^{r-2}s_{i}\choose
z-1}(z-1)_{(m-rl-t-h)\uparrow1}\mathscr {C}\bigl(z-1,j-r;\sigma,-(j-r)\sigma\bigr)
\\
&&\qquad {} \times\frac{1}{\sigma^{j-r}}\bigl((1-\sigma)_{(l-1)\uparrow
1}\bigr)^{2}2!(-1)^{2}
\mathscr{C}\Biggl(n-z+1-\sum_{i=1}^{r-2}s_{i},2;
\sigma -l,(j-r)\sigma-h\Biggr).
\end{eqnarray*}
The resulting expression has the same structure of the summand in
\eqref{eq_sum_r1}. This fact suggests the possibility of repeating
exactly the above arguments to each of the remaining nested sum over
the indexes $s_{r-2},\ldots,s_{1}$, respectively. In particular, after
a repeated application of these arguments we can write the $r$th
factorial moment of $R_{l,m}^{(n,j)}$ in \eqref{eq:starting3} as
%
\begin{eqnarray}
\label{eq:starting4} &&\E \bigl[ \bigl(R_{l,m}^{(n,j)}
\bigr)_{r\downarrow1} \bigr]\nonumber
\\
&&\quad
=\frac{\sigma^{j}}{\mathscr{C}(n,j;\sigma)}{m\choose l,\ldots,l,m-rl}r!\frac{(-(1-\sigma)_{l-1})^{r}}{\sigma^{j-r}}\sum
_{k=0}^{m}\frac{V_{n+m,j+k}}{V_{n,j}}
\frac{1}{\sigma^{k}}\nonumber
\\[-8pt]\\[-8pt]
&&\qquad {}
\times\sum_{t=k}^{m-rl}{m-rl
\choose t}\mathscr {C}(t,k;\sigma)\sum_{h=0}^{m-rl-t}{m-rl-t
\choose h}\bigl(-(j-r)\sigma \bigr)_{h\uparrow1}\nonumber
\\
&&\qquad {}
\times\sum_{z=j-r}^{n-r}{n\choose
z}(z)_{(m-rl-t-h)\uparrow
1}\mathscr{C}\bigl(z,j-r;\sigma,-(j-r)\sigma\bigr)\mathscr{C}
\bigl(n-z,r;\sigma -l,(j-r)\sigma-h\bigr).\nonumber
\end{eqnarray}
Finally, by applying \eqref{eq:gen_fact_coeff1} to expand $\mathscr
{C}(n-z,r;\sigma-l,(j-r)\sigma-h)$ in \eqref{eq:starting4}, we can
write \eqref{eq:starting4} as
\begin{eqnarray*}
&&\E \bigl[ \bigl(R_{l,m}^{(n,j)} \bigr)_{r\downarrow1} \bigr]
\\
&&\quad  =\frac{\sigma^{j}}{\mathscr{C}(n,j;\sigma)}{m\choose l,\ldots ,l,m-rl}r!\frac{(-(1-\sigma)_{l-1})^{r}}{\sigma^{j-r}}\sum
_{k=0}^{m}\frac{V_{n+m,j+k}}{V_{n,j}}
\frac{1}{\sigma^{k}}
\\
&&\qquad {} \times\sum_{s=r}^{n-(j-r)}{n\choose s}
\mathscr{C}(s,r;\sigma -l)\sum_{t=k}^{m-rl}{m-rl
\choose t}\\
&&\qquad {}\times\mathscr{C}(t,k;\sigma)\sum_{h=0}^{m-rl-t}{m-rl-t
\choose h}\bigl(-(j-r)\sigma\bigr)_{h\uparrow1}
\\
&&\qquad {} \times\sum_{z=j-r}^{n-s}{n-s\choose
z}(z)_{(m-rl-t-h)\uparrow
1}\bigl(-(j-r)\sigma+h\bigr)_{(n-s-z)\uparrow1}\\
&&\qquad {}\times\mathscr{C}
\bigl(z,j-r;\sigma ,-(j-r)\sigma\bigr)
\end{eqnarray*}
which leads to \eqref{teo2_eq2} by means on \eqref
{eq:gen_fact_coeff1} and some standard algebra involving factorial
numbers and noncentral generalized factorial coefficients. This
completes the second part of the proof.
\end{pf*}
\begin{pf*}{Proof of Proposition~\ref{prop2_1}} By combining the $r$th
factorial moment of $R_{l,m}^{(n,j,\mathbf{n})}$ in Theorem~\ref
{teo2} with $V_{n,j}$ displayed in \eqref{eq:defi_ewpit} one has
%
\begin{eqnarray}
\label{eq:start_freq_prop1} &&\E \bigl[ \bigl(R_{l,m}^{(n,j,\mathbf{n})}
\bigr)_{r\downarrow
1} \bigr]\nonumber
\\
&&\quad
=\frac{r!}{(\theta+n)_{m\uparrow1}}{m\choose l,\ldots ,l,m-rl}
\nonumber\\
&&\qquad {}
\times\sum_{\{c_{1},\ldots,c_{r}\}\in\mathcal
{C}_{j,r}}\prod
_{i=1}^{r}(n_{c_{i}}-\sigma)_{l\uparrow1}
\nonumber\\[-8pt]\\[-8pt]
&&\qquad {}
\times\sum_{k=0}^{m} \biggl(
\frac{\theta}{\sigma}+j \biggr)_{k\uparrow1}\mathscr{C}\Biggl(m-rl,k;\sigma,-n+
\sum_{i=1}^{r}n_{c_{i}}+(j-r)\sigma
\Biggr)\nonumber
\\
&&\quad
=\frac{r!}{(\theta+n)_{m\uparrow1}}{m\choose l,\ldots ,l,m-rl}
\nonumber\\
&&\qquad {}
\times\sum_{\{c_{1},\ldots,c_{r}\}\in\mathcal
{C}_{j,r}}\prod
_{i=1}^{r}(n_{c_{i}}-\sigma)_{l\uparrow1}
\Biggl(\theta +n-\sum_{i=1}^{r}n_{c_{i}}+
\sigma r\Biggr)_{(m-rl)\uparrow1},\nonumber
\end{eqnarray}
where the last identity follows equation 2.49 in \cite{Cha05}.
Accordingly, \eqref{eq:est1_freq_twopar} follows from \eqref
{eq:start_freq_prop1} by setting $r=1$. With regard to \eqref
{eq:est1_freq_twopar}, an inversion of the generating function for the
$r$th factorial moment in \eqref{eq:start_freq_prop1} leads to\vspace*{-1pt}
%
\begin{eqnarray}
\label{eq:start1_freq_prop1} &&\P \bigl[R_{l,m}^{(n,j,\mathbf{n})}=x \bigr]
\nonumber\\
&&\quad
=\frac{1}{(\theta+n)_{m\uparrow1}}\sum_{y\geq0}
\frac
{1}{x!}\frac{\ddr^{x}}{\ddr t^{x}}(t-1)^{x+y}\biggl\vert
_{t=0}{m\choose l,\ldots,l,m-(x+y)l}
\\
&&\qquad {}
\times\sum_{\{c_{1},\ldots,c_{x+y}\}\in\mathcal
{C}_{j,x+y}}\prod
_{i=1}^{x+y}(n_{c_{i}}-\sigma)_{l\uparrow1}
\Biggl(\theta +n-\sum_{i=1}^{x+y}n_{c_{i}}+
\sigma(x+y)\Biggr)_{(m-(x+y)l)\uparrow1},\nonumber
\end{eqnarray}
where\vspace*{-1pt}
%
\begin{eqnarray*}
\frac{\ddr^{x}}{\ddr t^{x} }(t-1)^{x+y}\biggl\vert _{t=0}=(-1)^{y}(x+y)_{x\downarrow1}.
\end{eqnarray*}
Then \eqref{eq:prob1_freq_twopar} follows from \eqref
{eq:start1_freq_prop1} by means of standard algebra involving factorial
numbers and binomial coefficients.
\end{pf*}
\begin{pf*}{Proof of Proposition~\ref{prop2_2}} A combination of the
$r$th factorial moment of $R_{l,m}^{(n,j)}$ in Theorem \ref{teo2}
with $V_{n,j}$ displayed in \eqref{eq:defi_ewpit} leads to\vspace*{-1pt}
%
\begin{eqnarray}
\label{eq:start_freq_prop2} &&\E \bigl[ \bigl(R_{l,m}^{(n,j)}
\bigr)_{r\downarrow1} \bigr]\nonumber
\\
&&\quad
=\frac{1}{\mathscr{C}(n,j;\sigma)(\theta+n)_{m\uparrow
1}}{m\choose l,\ldots,l,m-rl}r!\bigl(-\sigma(1-
\sigma)_{(l-1)\uparrow
1}\bigr)^{r}\nonumber
\\
&&\qquad {}
\times\sum_{s=r}^{n-(j-r)}{n\choose
s}\mathscr {C}(s,r;\sigma-l)\mathscr{C}(n-s,j-r;\sigma)
\nonumber\\
&&\qquad {}
\times\sum_{k=0}^{m} \biggl(
\frac{\theta}{\sigma}+j \biggr)_{k\uparrow1}\mathscr{C}\bigl(m-rl,k;
\sigma,-n+s+(j-r)\sigma\bigr)
\\
&&\quad
=\frac{1}{\mathscr{C}(n,j;\sigma)(\theta+n)_{m\uparrow
1}}{m\choose l,\ldots,l,m-rl}r!\bigl(-\sigma(1-
\sigma)_{(l-1)\uparrow
1}\bigr)^{r}\nonumber
\\
&&\qquad {}
\times\sum_{s=r}^{n-(j-r)}{n\choose
s}(\theta+n-s+\sigma r)_{(m-rl)\uparrow1}\nonumber\\
&&\qquad {}\mathscr{C}(s,r;\sigma-l)\mathscr
{C}(n-s,j-r;\sigma),\nonumber
\end{eqnarray}
where the last identity follows from equation 2.49 in \cite{Cha05}.
Accordingly, \eqref{eq:est2_freq_twopar} follows from \eqref
{eq:start_freq_prop2} by setting $r=1$. With regard to \eqref
{eq:prob2_freq_twopar}, an inversion of the generating function for the
$r$th factorial moment in \eqref{eq:start_freq_prop2} leads to
%
\begin{eqnarray}
\label{eq:start1_freq_prop2} &&\P \bigl[R_{l,m}^{(n,j)}=x \bigr]
\nonumber\\
&&\quad
=\frac{1}{\mathscr{C}(n,j;\sigma)(\theta+n)_{m\uparrow
1}}\sum_{y\geq0}
\frac{1}{x!}\frac{\ddr^{r}}{\ddr
t^{x}}(t-1)^{x+y}\biggl\vert
_{t=0}\nonumber
\\
&&\qquad {}
\times{m\choose l,\ldots,l,m-(x+y)l}\bigl(-\sigma(1-\sigma
)_{(l-1)\uparrow1}\bigr)^{(x+y)}
\\
&&\qquad {}
\times\sum_{s=x+y}^{n-(j-x-y)}{n\choose
s}\bigl(\theta +n-s+\sigma(x+y)\bigr)_{(m-(x+y)l)\uparrow1}\nonumber
\\
&&\qquad {}
\times\mathscr{C}(s,x+y;\sigma-l)\mathscr {C}\bigl(n-s,j-(x+y);\sigma
\bigr),\nonumber
\end{eqnarray}
where
%
\begin{eqnarray*}
\frac{\ddr^{x}}{\ddr t^{x} }(t-1)^{x+y}\biggl\vert _{t=0}=(-1)^{y}(x+y)_{x\downarrow1}.
\end{eqnarray*}
Then \eqref{eq:prob2_freq_twopar} follows from \eqref
{eq:start1_freq_prop2} by means of standard algebra involving factorial
numbers and binomial coefficients.
\end{pf*}
\subsection{Proofs of the results in Section \texorpdfstring{\protect\ref{subsec3}}{3.3}}

\begin{pf*}{Proof of Lemma~\ref{margin_lem}} By suitably marginalizing
the EPPF in \eqref{eq:gibbs_defi} one obtains the distribution of
$(K_{n},\mathbf{N}_{\tau})$, that is the main ingredient for
determining \eqref{eq:marginal_eppf}. Specifically, one has
%
\begin{eqnarray}
\label{eq:marginal} &&\P[K_{n}=j,\mathbf{N}_{\tau,n}=
\mathbf{n}_{\tau}]\nonumber
\\
&&\quad
=V_{n,j}\frac{(j-p)!}{j!}{n\choose n_{\tau_{1}},\ldots
,n_{\tau_{p}},n-\sum_{i=1}^{p}n_{\tau_{i}}}
\prod_{i=1}^{p}(1-\sigma )_{(n_{\tau_{i}}-1)\uparrow1}\nonumber
\\
&&\qquad {}
\times\frac{1}{(j-p)!}\sum_{(n_{\nu_{1}},\ldots,n_{\nu
_{j-p}})\in\mathcal{D}_{n-\sum_{i=1}^{p}n_{\tau_{i}},j-p}}{n-
\sum_{i=1}^{p}n_{\tau_{i}}\choose
n_{\nu_{1}},\ldots,n_{\nu
_{j-p}}}\nonumber\\[-8pt]\\[-8pt]
&&\qquad {}\times\prod_{i=1}^{(j-p)}(1-
\sigma)_{(n_{\nu_{i}}-1)\uparrow1}
\nonumber\\
&&\quad
=V_{n,j}\frac{(j-p)!}{j!}{n\choose n_{\tau_{1}},\ldots
,n_{\tau_{p}},n-\sum_{i=1}^{p}n_{\tau_{i}}}
\prod_{i=1}^{p}(1-\sigma )_{(n_{\tau_{i}}-1)\uparrow1}
\nonumber\\
&&\qquad {}
\times\frac{\mathscr{C}(n-\sum_{i=1}^{p}n_{\tau
_{i}},j-p;\sigma)}{\sigma^{j-p}}\nonumber
\end{eqnarray}
where the last identity is obtained by a direct application of equation
2.61 in \cite{Cha05}. The proof is completed by taking the ratio
between the distributions displayed in \eqref{eq:gibbs_defi} and
\eqref{eq:marginal}.
\end{pf*}
\end{appendix}

\section*{Acknowledgements} The authors are grateful to an Associate
Editor and a Referee for valuable remarks and suggestions that have
lead to a substantial improvement in the presentation. Stefano Favaro
is supported by the European Research Council (ERC) through StG
``N-BNP'' 306406.

%




\printhistory

\end{document}